\documentclass[12pt]{article}
\usepackage[dvips]{graphicx}
\usepackage{amsmath,amssymb,amsthm, color}
\def\frk{\mathfrak}               

\def\Phi{{\frk N}}
\def\opn#1#2{\def#1{\operatorname{#2}}} 
\opn\chara{char} 
\opn\length{\ell} 
\opn\pd{pd} 
\opn\rk{rk}
\opn\projdim{proj\,dim} 
\opn\injdim{inj\,dim} 
\opn\rank{rank}
\opn\depth{depth} 
\opn\grade{grade} 
\opn\height{height}
\opn\embdim{emb\,dim} 
\opn\codim{codim}

\opn\Tr{Tr} 
\opn\bigrank{big\,rank}
\opn\superheight{superheight}
\opn\lcm{lcm}
\opn\trdeg{tr\,deg}
\opn\reg{reg} 
\opn\lreg{lreg} 
\opn\ini{in} 
\opn\lpd{lpd}
\opn\size{size}
\opn\mult{mult}
\opn\dist{dist}
\opn\cone{cone}
\opn\lex{lex}
\opn\rev{rev}
\opn\div{div} \opn\Div{Div} \opn\cl{cl} \opn\Cl{Cl}
\opn\Spec{Spec} \opn\Supp{Supp} \opn\supp{supp} \opn\Sing{Sing}
\opn\Ass{Ass} \opn\Min{Min}
\opn\Ann{Ann} \opn\Rad{Rad} \opn\Soc{Soc}
\opn\Syz{Syz} \opn\Im{Im} \opn\Ker{Ker} \opn\Coker{Coker}
\opn\Am{Am} \opn\Hom{Hom} \opn\Tor{Tor} \opn\Ext{Ext}
\opn\End{End} \opn\Aut{Aut} \opn\id{id} \opn\ini{in}

\opn\nat{nat}
\opn\pff{pf}
\opn\Pf{Pf} \opn\GL{GL} \opn\SL{SL} \opn\mod{mod} \opn\ord{ord}
\opn\Gin{Gin}
\opn\Hilb{Hilb}\opn\adeg{adeg}\opn\std{std}\opn\ip{infpt}
\opn\Pol{Pol}
\opn\sat{sat}
\opn\Var{Var}
\opn\Gen{Gen}

\opn\aff{aff} \opn\con{conv} \opn\relint{relint} \opn\st{st}
\opn\lk{lk} \opn\cn{cn} \opn\core{core} \opn\vol{vol}
\opn\link{link} \opn\star{star}
\opn\gr{gr}

\def\pot#1#2{#1[\kern-0.28ex[#2]\kern-0.28ex]}

\opn\dirlim{\underrightarrow{\lim}}
\opn\inivlim{\underleftarrow{\lim}}
\def\Implies{\ifmmode\Longrightarrow \else
        \unskip${}\Longrightarrow{}$\ignorespaces\fi}
\def\implies{\ifmmode\Rightarrow \else
        \unskip${}\Rightarrow{}$\ignorespaces\fi}
\def\iff{\ifmmode\Longleftrightarrow \else
        \unskip${}\Longleftrightarrow{}$\ignorespaces\fi}

\let\:=\colon
\newtheorem{Theorem}{Theorem}[section]

\newtheorem{Corollary}[Theorem]{Corollary}

\theoremstyle{definition}

\newtheorem{Remark}[Theorem]{Remark}
\newtheorem{Remarks}[Theorem]{Remarks}
\newtheorem{Example}[Theorem]{Example}

\newtheorem{Definition}[Theorem]{Definition}

\let\epsilon\varepsilon
\let\phi=\varphi
\let\kappa=\varkappa
\def\qed{\ifhmode\textqed\fi
      \ifmmode\ifinner\quad\qedsymbol\else\dispqed\fi\fi}
\def\textqed{\unskip\nobreak\penalty50
       \hskip2em\hbox{}\nobreak\hfil\qedsymbol
       \parfillskip=0pt \finalhyphendemerits=0}
\def\dispqed{\rlap{\qquad\qedsymbol}}

\def\Bt{{\bf t}}
\def\Bx{{\bf x}}
\def\BX{{\bf X}}
\def\By{{\bf y}}
\def\Ba{{\bf a}}

\def\Bu{{\bf u}}
\def\Bv{{\bf v}}
\def\supp{{\rm supp}}

\def\Bt{{\bf t}}
\def\BT{{\bf T}}
\def\Bz{{\bf z}}
\def\Btheta{{\bf \theta}}

\def\Bpi{{\bf \pi}}
\def\Bone{{\bf 1}}
\def\Bzero{{\bf 0}}

\title{{\Large
An introduction to computational algebraic statistics
}}
\author{Satoshi Aoki}

\begin{document}
\maketitle

\begin{abstract}
In this paper, we introduce the fundamental notion of a Markov basis,
 which is one of the first connections between commutative algebra and
 statistics. The notion of a Markov basis is first introduced by Diaconis
 and Sturmfels (\cite{Diaconis-Sturmfels-1998}) for conditional testing
 problems on contingency tables by Markov chain Monte Carlo methods. In
 this method, we make use of a connected Markov
 chain over the given conditional sample space to 
 estimate the $p$ values numerically for various conditional tests. 
A Markov basis plays an importance role in this arguments, because it 
guarantees the connectivity of the chain, which is needed for
 unbiasedness of the estimate, for
 arbitrary conditional sample space. As another important point, 
a Markov basis is characterized as generators of
 the well-specified toric ideals of polynomial rings. 
This connection between commutative algebra and statistics is the main
 result of \cite{Diaconis-Sturmfels-1998}. After this first paper, 
a Markov basis is studied intensively by 
many researchers both in commutative algebra and statistics, which
 yields an attractive field called {\it computational algebraic
 statistics}. 
In this paper, we give a review of the Markov chain Monte
 Carlo methods for contingency tables and Markov bases, with some 
fundamental examples. We also give some computational examples by
 algebraic software Macaulay2
 (\cite{Macaulay2}) and statistical software \verb|R|.
Readers can also find theoretical details of the problems considered in
 this paper 
and various results on the structure and examples of Markov
 bases in \cite{Aoki-Hara-Takemura-2012}.
\end{abstract}

\section{Conditional tests for contingency tables}
A contingency table is a cross-classified table of frequencies. For
example, suppose $40$ students in some class took examinations of two
subjects, Algebra and Statistics. Suppose that both scores are classified to
one of the categories, \{Excellent, Good, Fair\}, and 
are summarized in Table \ref{tbl:3x3-example}. 
\begin{table}[h]
\begin{center}
\caption{Scores of Algebra and Statistics for $40$ students (imaginary data)}
\label{tbl:3x3-example}
\begin{tabular}{l|r|r|r|r}
\multicolumn{1}{c}{Alg$\backslash$Stat} & \multicolumn{1}{c}{Excellent} & 
\multicolumn{1}{c}{\ \ Good\ \ \ } & \multicolumn{1}{c}{\ \ Fair\ \ \ \
 } & \multicolumn{1}{c}{Total}\\ \cline{2-4} 
Excellent & $11$ & $5$ & $2$ & $18$\\ \cline{2-4}
Good      & $4$  & $9$ & $1$ & $14$\\ \cline{2-4}
Fair      & $2$  & $3$ & $3$ & $8$\\ \cline{2-4}
\multicolumn{1}{l}{Total} & \multicolumn{1}{r}{$17$} &
	 \multicolumn{1}{r}{$17$} &
	     \multicolumn{1}{r}{$6$} & \multicolumn{1}{r}{$40$}
\end{tabular}
\end{center}
\end{table}
This is a typical example of two-way contingency tables. 
Since this table has $3$ rows and $3$ columns, this is
called a $3\times 3$ contingency table. 
The two subjects, Algebra and Statistics, are called {\it factors} of the
table, 
and the outcomes (i.e., scores) of each factor, \{Excellent, Good, Fair\}, are
called {\it levels} of each factor. The {\it cells} of the $I \times J$
contingency table is
the $IJ$ possible combinations of outcomes. 
Three-way, four-way or higher dimensional
contingency tables are defined similarly. 
For example, adding to the data of Table 
\ref{tbl:3x3-example}, if
the scores of another subject (Geometry, for example) are also given, we
have a three-way contingency table. An  
$I_1\times \cdots \times I_m$ ($m$-way)contingency table has 
$\prod_{i=1}^mI_i$ cells, where $I_i$ is the number of levels for the
$i$th factor, $i = 1,\ldots,m$. 
In statistical data analysis, the development of methods for analyzing
contingency tables began in the 1960s. We refer to \cite{Agresti-2013}
for standard textbook in this field.

We begin with simple $I\times J$ cases, and will consider
generalizations to $m$-way cases afterward.
In statistical inference, we consider underlying random variables and
statistical models for observed data such as Table
\ref{tbl:3x3-example}, and treat the observed data as one realization of
the random variables. In the case of Table \ref{tbl:3x3-example}, it is
natural to deal with the two-dimensional discrete random variables 
\begin{equation}
(V_1, W_1), (V_2, W_2),\ldots,(V_n,W_n), 
\label{eqn:2-dim-rv}
\end{equation}
where $n$ is the {\it sample size}, ($n=40$ for Table
\ref{tbl:3x3-example}) and 
$(V_k,W_k)$ is the couple of 
scores obtained by the $k$th student. The random couples $(V_k, W_k)$
for $k = 1,\ldots,n$ are drawn
independently from
the same distribution
\[
 P(V_k = i,\ W_k = j) = \theta_{ij},\ \  i \in [I],\ j \in [J],\ k
 \in [n].
\]
Here we use a notation $[r] = \{1,2,\ldots,r\}$ for $r \in
\mathbb{Z}_{\geq 0}$, where $\mathbb{Z}_{\geq 0}$ is the set of
nonnegative integers. Note that we use appropriate coding such as $1$:\
Excellent, $2$:\ Good, $3$:\ Fair. The probability $\Btheta =
(\theta_{ij})$ satisfies the condition
\[
 \sum_{i=1}^{I}\sum_{j=1}^{J}\theta_{ij} = 1,
\]
and is 
called a {\it parameter}. The parameter space
\[
 \Delta_{IJ-1} = \left\{(\theta_{11},\ldots,\theta_{IJ}) \in
      \mathbb{R}^{IJ}_{\geq 0}:\ \sum_{i=1}^{I}\sum_{j=1}^{J}
     \theta_{ij} = 1\right\}
\]
is called an $IJ-1$ dimensional {\it probability simplex}.

To consider the data in the form of a contingency table, 
we also summarize the underlying random variable
(\ref{eqn:2-dim-rv}) to the form of the contingency tables as
\[
 X_{ij} = \sum_{k = 1}^n\Bone(V_k = i, W_k = j),
\]
for $i \in [I], j \in [J]$, where $\Bone(\cdot)$ is the indicator
function. By this aggregation from the raw scores to the contingency
table, we neglect the order of
observations in (\ref{eqn:2-dim-rv}), that is considered to have no
information for estimating the parameter $\Btheta$. 
Then the data $\Bx = (x_{ij}) \in \mathbb{Z}^{IJ}_{\geq 0}$ is treated as a
realization of $\BX = (X_{ij})$. 
The distribution of $\BX$ is a {\it multinomial distribution} given by
\begin{equation}
 p(\Bx) = P(\BX = \Bx) =
      \frac{n!}{\displaystyle\prod_{i=1}^I\prod_{j=1}^Jx_{ij}!}
\prod_{i=1}^I\prod_{j=1}^J\theta_{ij}^{x_{ij}},\ \ 
 \sum_{i=1}^I\sum_{j=1}^Jx_{ij} = n.
\label{eqn:multinom-dist}
\end{equation}
We see that the multinomial distribution (\ref{eqn:multinom-dist}) is
derived from the joint probability function for $n$ individuals
under the
assumption that each outcome is obtained independently.

By summarizing the data in the form of contingency tables for fixed
sample size $n$, 
the degree of freedom of the observed frequency $\Bx$
 becomes $IJ - 1$. which coincides the degree of freedom of
the parameter $\Btheta \in \Delta_{IJ-1}$. Here, we use 
``degree of freedom'' as the number of elements that are free to vary, 
that is a well-used terminology in statistical fields.  
We can see 
the probability simplex $\Delta_{IJ-1}$ as an example of 
statistical models, called a {\it saturated model}. 
Statistical model is called saturated if the degree of freedom of the
      parameter equals to the degree of freedom of data.

The saturated model is also characterized as the statistical model having the
 parameter with the largest degree of freedom.  
In this sense, the saturated model is the most complex statistical model. 
In other words, the saturated model is the statistical model that fits
 the observed data perfectly, i.e., fits the data {\it without error}. 
In fact, the
 parameter $\Btheta$ in the saturated model $\Delta_{IJ-1}$ is estimated
 from the data as
\begin{equation}
 \hat{\theta}_{ij} = \frac{x_{ij}}{n},\ i \in [I],\ j \in [J],
\label{eqn:empirical-2-way}
\end{equation}
that is also called an empirical probability of data. Because we
assume that 
the data $\Bx$ is obtained from some probability function such as
multinomial distribution (\ref{eqn:multinom-dist}) {\it with some randomness},
we want to consider more simple statistical model, i.e., a subset of the
saturated model,  
${\cal M} \subset \Delta_{IJ - 1}$.

In the two-way contingency tables, a natural, representative statistical
model is an {\it independence model}. 
\begin{Definition}
The independence model for $I\times J$ contingency tables is the set
\begin{equation}
 {\cal M}_{indp} = \{\Btheta \in \Delta_{IJ-1}\ :\ \theta_{ij} =
 \theta_{i+}\theta_{+j}, \forall i, \forall j\},
\label{eqn:indep-model}
\end{equation}
where
\[
 \theta_{i+} = \sum_{j = 1}^J\theta_{ij},\ 
 \theta_{+j} = \sum_{i = 1}^I\theta_{ij}\ \ \mbox{for}\ i \in [I], j \in [J].
\]
\end{Definition}

\begin{Remarks}
Here we consider that only the sample size $n$ is
 fixed. However, several different situations can be considered for
 $I\times J$ contingency tables. The situation that we consider here is
 called a multinomial sampling scheme. For other sampling schemes such
 as Poisson, binomial and so on, see 
Chapter 2 of \cite{Agresti-2013} or Chapter 4 of \cite{dojo-en}. 
Accordingly, 
the corresponding independence model ${\cal M}_{indp}$ is called in
 different way for
 other sampling schemes. For example, it is called a {\it common
 proportions model} for (product of) binomial sampling scheme where the
 row sums are fixed, and {\it
 main effect model} for Poisson sampling scheme where no marginal is
 fixed. Though there are also a little
 differences between the descriptions of these models, we can treat
 these models almost in the same way by considering the {\it conditional
 probability function}, which we consider afterward. Therefore we
 restrict our arguments to the multinomial sampling scheme in this paper.
\end{Remarks}

There are several equivalent descriptions for the independence model
${\cal M}_{indp}$. The most common 
parametric description in statistical textbooks is 
\begin{equation}
 {\cal M}_{indp} = \left\{\Btheta \in \Delta_{IJ-1}:\ 
\theta_{ij} =
      \alpha_{i}\beta_{j}\ {\rm for\ some}\ (\alpha_i),\ (\beta_j)
\right\}.
\label{eqn:parametric-description}
\end{equation}
For other equivalent parametric descriptions or implicit
descriptions, see Section 1 of \cite{Pacher-Sturmfels}, for example.
%

The meaning of ${\cal M}_{indp}$ in Table \ref{tbl:3x3-example} is 
as follows. If ${\cal M}_{indp}$ is true, there are no relations between
the scores of two subjects. Then we can imagine that the scores of two 
subjects follow the marginal probability functions for each score
respectively, and are independent, and the discrepancy we observed in Table 
\ref{tbl:3x3-example} is obtained ``by chance''. 
However, it is natural to imagine some
structure between the two scores such as ``there is a tendency that the
students having better scores in Algebra are likely to have better scores
in Statistics'', because these subjects are in the same mathematical
category. In fact, we see relatively large frequencies $11$ and $9$ in
the diagonals of Table \ref{tbl:3x3-example}, which seem to indicate a
positive correlation. 
Therefore one of natural questions for Table \ref{tbl:3x3-example} is 
``Is there some tendency between the two scores that breaks independence?''.
To answer this question, 
we evaluate the fitting of ${\cal M}_{indp}$
by {\it hypothetical testing}.

The hypothetical testing problem that we consider in this paper
is as follows.
\[
 \mbox{H}_0:\ \Btheta \in {\cal M}_{indp}\ \ \mbox{v.s.}\ \ 
\mbox{H}_1:\ \Btheta \in \Delta_{IJ-1}\setminus {\cal M}_{indp}.
\]
Here we call $\mbox{H}_0$ a {\it null hypothesis} and 
$\mbox{H}_1$ an {\it alternative hypothesis}. The terms {\it
null model} and {\it alternative model} are also used. 
The hypothetical testing
in the above form, i.e., a null model is a subset of a saturated
model, ${\cal M} \subset \Delta_{IJ-1}$, and the alternative model is
the complementary set of ${\cal M}$ into the saturated model, 
is called a {\it goodness-of-fit test} of
model ${\cal M}$. The testing procedures are composed of steps such as
choosing a test statistics, choosing a significance level, and
calculating the $p$ value. We see these steps in order.

\paragraph*{Choosing a test statistic.}\ First we have to choose a test
statistic to use. In general, the term {\it statistic} means a
function of the
random variable $\BX = (X_{ij})$. For example, $(X_{i+})$ and $(X_{+j})$
given by
\[
 X_{i+} = \sum_{j = 1}^J X_{ij},\ 
 X_{+j} = \sum_{i = 1}^I X_{ij}\ \mbox{for}\ i \in [I], j \in [J]
\]
are examples of statistics called the row sums and the column sums,
respectively. Other examples of statistics are the row mean $\bar{X}_{i+} =
X_{i+}/J$ and the column mean $\bar{X}_{+j}/I$ for $i \in [I], j \in
[J]$. To perform the hypothetical testing, we first select an appropriate
statistic, called a {\it test statistic}, to measure the discrepancy
of the observed data from the null model. 
One of the common test statistic for the goodness-of-fit test is
a {\it Pearson goodness-of-fit $\chi^2$} given by
\[
 \chi^2(\BX) = \sum_{i = 1}^I\sum_{j = 1}^J\frac{(X_{ij} -
 \hat{m}_{ij})^2}{\hat{m}_{ij}},
\]
where $\hat{m}_{ij}$ is the {\it fitted value} of $X_{ij}$ under $\mbox{H}_0$, 
i.e., an estimator of $E(X_{ij}) = m_{ij} = n\theta_{ij}$, 
given by
\begin{equation}
 \hat{m}_{ij} = n\hat{\theta}_{ij} = \frac{x_{i+}x_{+j}}{n}.
\label{eqn:fitted-value-m}
\end{equation}
Here we use the {\it maximum likelihood estimate} of the parameter
under the null model,
$\hat{\Btheta} = (\hat{\theta}_{ij})$, given by
\begin{equation}
 \hat{\theta}_{ij} = \frac{x_{i+}x_{+j}}{n^2},
\label{eqn:MLE-under-indp-thetaij}
\end{equation}
that is obtained by maximizing the log-likelihood
\[
 \mbox{Const} + \sum_{i = 1}^I\sum_{j = 1}^Jx_{ij}\log\theta_{ij}
\]
under the constraint $\Btheta \in {\cal M}_{indp}$. The meaning of this
estimate is also clear in a parametric description
(\ref{eqn:parametric-description}) since the maximum likelihood estimates
of $(\alpha_i), (\beta_j)$ are given by
\[
 \hat{\alpha}_i = \frac{x_{i+}}{n},\ \ 
 \hat{\beta}_j = \frac{x_{+j}}{n},
\]
respectively. The fitted value for Table \ref{tbl:3x3-example} under
${\cal M}_{indp}$ is given in Table \ref{tbl:mle-3x3-example}.
\begin{table}[h]
\begin{center}
\caption{The fitted value under ${\cal M}_{indp}$ for Table
 \ref{tbl:3x3-example}}
\label{tbl:mle-3x3-example}
\begin{tabular}{l|r|r|r|r}
\multicolumn{1}{c}{Alg$\backslash$Stat} & \multicolumn{1}{c}{Excellent} & 
\multicolumn{1}{c}{\ \ Good\ \ \ } & \multicolumn{1}{c}{\ \ Fair\ \ \ \
 } & \multicolumn{1}{c}{Total}\\ \cline{2-4} 
Excellent & $7.65$ & $7.65$ & $2.70$ & $18$\\ \cline{2-4}
Good      & $5.95$  & $5.95$ & $2.10$ & $14$\\ \cline{2-4}
Fair      & $3.40$  & $3.40$ & $1.20$ & $8$\\ \cline{2-4}
\multicolumn{1}{l}{Total} & \multicolumn{1}{r}{$17$} &
	 \multicolumn{1}{r}{$17$} &
	     \multicolumn{1}{r}{$6$} & \multicolumn{1}{r}{$40$}
\end{tabular}
\end{center}
\end{table}

There are various test statistics other than the Pearson goodness-of-fit
$\chi^2$ that can be used in our problem. Another representative is the
(twice log) likelihood ratio given by
\begin{equation}
 2\sum_{i = 1}^I\sum_{j = 1}^JX_{ij}\log\frac{X_{ij}}{\hat{m}_{ij}},
\label{eqn:likelihood-ratio}
\end{equation}
where $\hat{m}_{ij}$ is given by (\ref{eqn:fitted-value-m}). In general, test
statistic should be selected by considering their {\it power}, i.e.,
the probability that the null hypothesis is rejected if the alternative
hypothesis is true. See textbooks such as \cite{Lehmann-TSH}
for the theory of the hypothetical testing, 
the optimality of the test statistics, examples and
the guidelines for choosing test statistics for various problems. 

\paragraph*{Choosing a significance level.}\ Once we choose a test
statistic to use, as the Pearson goodness-of-fit $\chi^2$ for example,
the hypothetical testing procedure is written by
\[
 \chi^2(\Bx^o) \geq c_{\alpha}\ \Rightarrow\ \mbox{Reject}\ \mbox{H}_0,
\]
where $\Bx^o$ is the observed data, and $c_{\alpha}$ is the {\it
critical point at the significance level $\alpha$} satisfying
\begin{equation}
 P(\chi^2(X) \geq c_{\alpha}\ |\ \mbox{H}_0) \leq \alpha.
\label{eqn:prob-type-1-error}
\end{equation}
The probability of the left hand side of (\ref{eqn:prob-type-1-error})
is called a {\it type I error}. 
Equivalently, we define the {\it $p$-value} by
\begin{equation}
 p = P(\chi^2(X) \geq \chi^2(\Bx^o)\ |\ \mbox{H}_0),
\label{eqn:p-value}
\end{equation}
then the testing procedure is written by
\[
 p \leq \alpha\ \Rightarrow\ \mbox{Reject}\ \mbox{H}_0.
\]
The meaning of the $p$-value for the data $\Bx^o$ is the conditional
probability that ``more or equally discrepant results are obtained than
the observed data if the null hypothesis is true''. 
Therefore, if $p$-value is significantly small, we conclude that null
hypothesis is unrealistic, because it is doubtful that 
such an extreme result $\Bx^o$ is obtained. This is the idea of the
statistical hypothetical testing. 
In this process, the significance level $\alpha$ plays a threshold to
decide the $p$-value is ``significantly small'' to reject the null
hypothesis.  In statistical and scientific literature, it is common to
choose $\alpha = 0.05$ or $\alpha = 0.01$. Readers can find various
topics on $p$-value in \cite{Vickers}.

\paragraph*{Calculating the $p$-value.}\ Once we choose a test
statistic and a significance level, all we have to do is to calculate the
$p$-value given in (\ref{eqn:p-value}) for observed data $\Bx^o$. 
The observed value of the Pearson goodness-of-fit $\chi^2$ 
for Table \ref{tbl:3x3-example} is
\[
\chi^2(\Bx^o) = \sum_{i=1}^3\sum_{j=1}^3\frac{(x_{ij}^o -
\hat{m}_{ij})^2}{\hat{m}_{ij}} = \frac{(11-7.65)^2}{7.65} + \cdots +
\frac{(3-1.20)^2}{1.20} = 8.6687,
\]
therefore the $p$-value for our $\Bx^o$ is
\[
 p = P(\chi^2(X) \geq 8.6687\ |\ \mbox{H}_0).
\]
This probability is evaluated based on the {\it probability function of
the test statistic $\chi^2(\BX)$ under H$_0$}, which we call a {\it
null distribution} hereafter. 
Unfortunately, the null distribution 
depends on the unknown parameter $\Btheta \in {\cal M}_{indp}$ and
the $p$-values cannot be calculated in most cases in principle. 
One naive idea to evaluate the $p$-values for such cases 
is to calculate its supremum in
${\cal M}_{indp}$ and perform the test as the form
\begin{equation}
 \sup_{\Btheta \in {\cal M}_{indp}} P(\chi^2(X) \geq \chi^2(\Bx^o)\ |\
 \mbox{H}_0) \leq \alpha\ \Rightarrow\ \mbox{Reject}\ \mbox{H}_0.
\label{eqn:sup}
\end{equation}
However, this idea is hard to implement in general, i.e., 
it is usually difficult to 
evaluate the left-hand side of
(\ref{eqn:sup}) or to seek tests that are powerful under (\ref{eqn:sup}). 
Then, what should we do? We consider the following three strategies for
calculating $p$-values in this paper. 
\begin{itemize}
\item[(a)] Using the asymptotic distribution of the test statistic.
\item[(b)] Exact calculation based on the conditional distribution.
\item[(c)] Estimate the $p$-value by the Monte Carlo method.
\end{itemize}
The aim of this paper is to introduce strategy (c). We will consider
each strategy in order.

\paragraph*{(a) Using the asymptotic distribution of the test statistic.}
In applications, it is common to rely on various asymptotic theories for
the test statistics. As
for the Pearson goodness-of-fit test $\chi^2$ test, the following result
is known.
\begin{Theorem}
Under the null model ${\cal M}_{indp}$, 
 the Pearson goodness-of-fit $\chi^2(X)$ asymptotically follows the $\chi^2$
 distribution with $(I-1)(J-1)$ degree of freedom, i.e.,
\[
 \lim_{n \rightarrow \infty}P(\chi^2(X) \geq u) = P(V \geq u)\ \
 \mbox{for}\ u>0,
\]
where $V \sim \chi^2_{(I-1)(J-1)}$, i.e., $V$ is distributed to the
 $\chi^2$ distribution with $(I-1)(J-1)$ degree of freedom.   
\label{thm:asymptotic-chi^2}
\end{Theorem}
This theorem is shown as a consequence of the central limit theorem. In
addition, the same asymptotic distribution is given when we consider
the conditional limit, i.e., consider $n \rightarrow \infty$ under the
condition that $X_{i+}/n \rightarrow a_{i}$ and $X_{+j}/n \rightarrow
b_j$ for $i \in [I], j \in [J]$ for some fixed $0 < a_i, b_j < 1$. 
See \cite{Cornfield} or \cite{Plackett} for detail.  
Anyway, these asymptotic properties are the reason why we call this test
as Pearson
goodness-of-fit ``$\chi^2$ test''. 
Similarly to the Pearson goodness-of-fit $\chi^2$, there are several
test statistics that have the $\chi^2$ distribution as the asymptotic
distribution. An important example is the likelihood ratio test
statistic, which is given in (\ref{eqn:likelihood-ratio}) for our
setting. Moreover, several asymptotic good properties of likelihood
ratio test statistics are known. See \cite{Lehmann-TSH} for
details. Note also that our methods, Markov chain Monte Carlo methods,
can be applicable for arbitrary type of test statistics, though we only
consider the Pearson goodness-of-fit $\chi^2$ in this paper. 

\begin{figure}[htbp]
\begin{center}
\includegraphics[width=70mm]{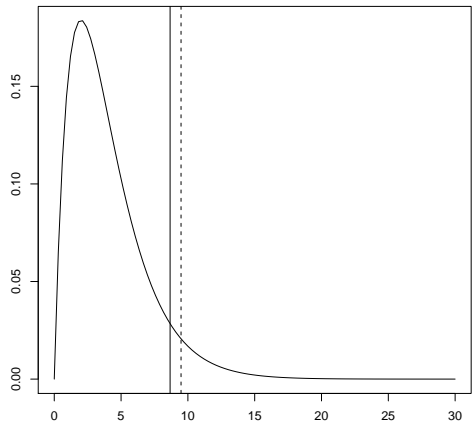}
\caption{$\chi^2$ distribution with degree of freedom $4$. The vertical
solid line indicates the observed value $\chi^2(\Bx^o) = 8.6687$, and the
 dotted line
 indicates the critical point for the significance level $\alpha =
 0.05$, $\chi^2_{4,0.05} = 9.488$.}
\label{fig:chi2-deg4}
\end{center}
\end{figure}

Following Theorem \ref{thm:asymptotic-chi^2}, it is easy to evaluate the
asymptotic $p$-value of the 
Pearson goodness-of-fit $\chi^2$ test. For our data, the observed value
of test statistic, $\chi^2(\Bx^o) = 8.6687$, is less than the upper
$5$ percent point of the $\chi^2$ distribution with $4$ degrees of
freedom, $\chi^2_{4, 0.05} = 9.488$. Therefore, for the significance level
$\alpha = 0.05$, we cannot reject the null hypothesis H$_0$, i.e., 
we cannot say that ``the fitting of the model ${\cal M}_{indp}$ to Table
\ref{tbl:3x3-example} is poor''. Equivalently, the asymptotic $p$-value
is calculated as the upper probability of $\chi^2_4$, which is 
$0.0699$ and is greater than $\alpha = 0.05$. Figure \ref{fig:chi2-deg4} 
presents the probability density function of the $\chi^2_4$ distribution.
The above results can be obtained numerically 
by the following codes of the statistical software \verb|R|.
\begin{verbatim}
> x <- matrix(c(11,5,2,4,9,1,2,3,3), byrow=T, ncol=3, nrow=3)
> x
     [,1] [,2] [,3]
[1,]   11    5    2
[2,]    4    9    1
[3,]    2    3    3
> chisq.test(x)

	Pearson's Chi-squared test

data:  x 
X-squared = 8.6687, df = 4, p-value = 0.06994
> pchisq(8.6687,4, lower.tail=F)               
[1] 0.06993543
> qchisq(0.05,4,lower.tail=F)                  # critical point
[1] 9.487729
\end{verbatim}

As we see above, using asymptotic null distribution is easy way to
evaluate $p$-values, and one of the most common approaches in
applications. One of the disadvantages of strategy (a) is that there
might not be a good fit with the asymptotic distribution. In fact,
because sample size is only $n=40$ for Table \ref{tbl:3x3-example}, it
is doubtful that we can apply the asymptotic result of $n \rightarrow
\infty$. Besides, it is well known that there is cases that 
the fitting of the asymptotic distributions are poor for data with
relatively large sample sizes. One such case is sparse data case,
another one is unbalanced case. See \cite{Haberman} for these topics. 

\paragraph*{(b) Exact calculation based on the conditional distribution.}
If we want to avoid asymptotic approaches as strategy (a),
an alternative 
choice is to calculate $p$-values {\it exactly}. For the cases that the null
distribution of the test statistics depend on the unknown parameters, 
we can formulate the exact methods based on the {\it conditional
probability functions} for fixed {\it minimal sufficient statistics}
under the null model ${\cal M}_{indp}$. The key notion here is the
minimal sufficient statistics.
\begin{Definition}
Let $\BX$ be a discrete random variable with the probability function
 $p(\Bx)$ with the parameter $\Btheta$.  
The statistic $\BT(\BX)$, i.e., a vector or a scalar function of $\BX$,
 is called sufficient for $\Btheta$ if the
 conditional probability function of $\BX$ for a given $\BT$,
\begin{equation}
 p(\Bx\ |\ \Bt) = P(\BX = \Bx\ |\ \BT(\BX) = \Bt),
\label{eqn:def-sufficient-statistics}
\end{equation}
does not depend on $\Btheta$. The sufficient statistic $\BT(\BX)$ is
 minimal if there is no other sufficient statistics that is a 
function of $\BT(\BX)$.  
\end{Definition}
The meaning of the minimal sufficient statistic is explained as
follows. If we know the value of $\BT$, then knowing $\BX$ provides no
further information about the parameter $\Btheta$. Therefore for the
parameter estimation or hypothetical testing, it is sufficient to
consider the methods based on the minimal sufficient statistic. 
The minimal sufficient statistics for our two-way problem is as follows.
\begin{itemize}
\item Under the saturated model $\Btheta \in {\Delta}_{IJ-1}$, a
      minimal sufficient statistic is the contingency table $\BX$. 
Adding the additional information such as the scores of the $k$th
      student, $(V_k, W_k)$, in (\ref{eqn:2-dim-rv}) gives us no
      additional information on the estimation of $\Btheta$. Indeed, 
      under the saturated model, the maximum likelihood estimate of the
      parameter is the empirical probability (\ref{eqn:empirical-2-way}),
      that is a function of the minimal sufficient statistic. 
\item Under the independence model ${\cal M}_{indp}$, a minimal
      sufficient statistic is the row sums $\{X_{i+},\ i \in [I]\}$ and
      the column
      sums $\{X_{+j},\ j \in [J]\}$, as we see below. Indeed,  
      we have already seen that the maximum likelihood estimate of
      the parameter under the independence model is
      (\ref{eqn:MLE-under-indp-thetaij}), that is a function of the row
      sums and column sums.
Note that $\BX$ itself is also the sufficient statistic under the
      independence model, but is not minimal. 
\end{itemize}

To see that a given statistic $\BT(\BX)$ is sufficient for a parameter
$\Btheta$, a useful way is to rely on the following theorem.
\begin{Theorem}
$\BT(\BX)$ is a sufficient statistic for $\Btheta$ if and only if the
 probability function of $\BX$ is factored as
\begin{equation}
 p(\Bx; \Btheta) = h(\Bx)g(T(\Bx); \Btheta),
\label{eqn:factorization-thm}
\end{equation}
where $g(\cdot)$ is a function that depends on the parameter $\Btheta$
 and $h(\cdot)$ is a function that does not.
\end{Theorem}
For the case of discrete probability function, this theorem, called a
{\it factorization theorem}, is easily (i.e., without measure theories) proved
from the definition of the sufficient statistic. 
Generally, to obtain such a factorization is easier than to compute
explicitly the conditional distribution
(\ref{eqn:def-sufficient-statistics}). 
For example, under the parametric description $\theta_{ij} =
\alpha_i\beta_j$, the probability function of the multinomial distribution 
(\ref{eqn:multinom-dist}) is written as
\[
 p(\Bx; \Btheta) = 
\frac{n!}{\prod\prod x_{ij}!}\left(\prod_i
 \alpha_i^{x_{i+}}\right)\left(\prod_j\beta_j^{x_{+j}}\right) 
\]
and we see that $T(\BX) = (\{X_{i+}\}, \{X_{+j}\})$ is a sufficient
statistic for the parameter $\Btheta \in {\cal M}_{indp}$. 

Here, for later generalization, we introduce a {\it configuration
matrix} $A$ and express a minimal sufficient statistic by $A$ as
follows. Let the number of the cells of the contingency table $\BX$ be
$p$ and treat $\BX$ as a $p$-dimensional column vector. 
Let $T(\BX)$ be a $d$-dimensional sufficient statistic for the
parameter $\Btheta \in {\cal M}$. 
For example of the independence model ${\cal M}_{indp}$ for $I\times J$
contingency tables, we have 
$\nu = IJ$, $\BX = (X_{11}, X_{12}, \ldots, X_{IJ})'$ and 
\[
\BT(\BX) = (X_{1+}, \ldots, X_{I+}, X_{+1}, \ldots, X_{+J})' 
\]
and $d = I + J$. Then we see that $T(\BX)$ is written as
\begin{equation}
 T(\BX) = A\BX
\label{eqn:t=Ax}
\end{equation}
for $d\times \nu$ integer matrix $A$. For the $3\times 3$ contingency
tables, $A$ is written as follows:
\begin{equation}
 A = \left(\begin{array}{ccccccccc}
1 & 1 & 1 & 0 & 0 & 0 & 0 & 0 & 0\\
0 & 0 & 0 & 1 & 1 & 1 & 0 & 0 & 0\\
0 & 0 & 0 & 0 & 0 & 0 & 1 & 1 & 1\\
1 & 0 & 0 & 1 & 0 & 0 & 1 & 0 & 0\\
0 & 1 & 0 & 0 & 1 & 0 & 0 & 1 & 0\\
0 & 0 & 1 & 0 & 0 & 1 & 0 & 0 & 1
\end{array}\right).
\label{eqn:A-3x3-indp}
\end{equation}

Following the sufficiency of $T(\BX) = A\BX$, the conditional
probability function for given $\BT = \Bt$ does not depend on the
parameter. For the case of the independence model ${\cal M}_{indp}$ for
two-way contingency tables, it is 
\begin{equation}
\begin{array}{l}
 h(\Bx)  = P(X = \Bx\ |\ A\BX = \Bt, \ \mbox{H}_0) 
 = 
\frac{\left(\displaystyle\prod_{i}x_{i+}!\right)
\left(\displaystyle\prod_{j}x_{+j}!\right)}
{n!\displaystyle\prod_{i,j}x_{ij}!},\ \ \ \Bx \in {\cal F}_{\Bt},
\label{eqn:hyper-geo}
\end{array}\end{equation}
where
\[
 {\cal F}_{\Bt} = \{\Bx \in \mathbb{Z}_{\geq 0}^{\nu}\ :\ A\Bx = \Bt
\}
\]
is the conditional sample space, which is called a {\it $\Bt$-fiber}
in the arguments of Markov bases. The conditional probability function
$h(\Bx)$ is called a {\it hypergeometric distribution}. Using this
conditional probability, the conditional $p$-value can be defined by
\begin{equation}
p = E_{\mbox{H}_0}(g(\BX)\ |\ A\BX = A\Bx^o) =  \sum_{\Bx \in {\cal
 F}_{A\Bx^o}}g(\Bx)h(\Bx)
\label{eqn:exact-p-value}
\end{equation}
for the observed table $\Bx^o$, where $g(\Bx)$ is the test function
\[
 g(\Bx) = \left\{\begin{array}{ll}
1, & \chi^2(\Bx) \geq \chi^2(\Bx^o),\\
0, & \mbox{otherwise}.
\end{array}
\right.
\]
Now calculate the conditional $p$-value exactly for Table
\ref{tbl:3x3-example}. 
For the observed table $\Bx^o$, i.e., Table
\ref{tbl:3x3-example}, we consider the independence model ${\cal
M}_{indp}$. The configuration matrix $A$ for ${\cal M}_{indp}$ is given
in (\ref{eqn:A-3x3-indp}). 
The $\Bt$-fiber including $\Bx^o$, i.e.,
$A\Bx^o$-fiber, is the set of all contingency tables that have the same
value of the row sums and the column sums to $\Bx^o$, 
\[
{\cal F}_{A\Bx^o} = \left\{
\Bx\in \mathbb{Z}_{\geq 0}^9 :\  \begin{array}{|c|c|c|r}\cline{1-3}
x_{11} & x_{12}& x_{13}& 18\\ \cline{1-3}
x_{21} & x_{22}& x_{23}& 14\\ \cline{1-3}
x_{31} & x_{32}& x_{33}& 8\\ \cline{1-3}
\multicolumn{1}{c}{17} & \multicolumn{1}{c}{17} & \multicolumn{1}{c}{6} &
\multicolumn{1}{c}{40}
\end{array}
\right\}\ .
\]
There are $2366$ elements in this ${\cal F}_{A\Bx^o}$. 
For each $2366$ elements in ${\cal F}_{A\Bx^o}$, the conditional
probability is given by
\[
 h(\Bx) = 
\frac{\left(\displaystyle 18!14!8!\right)
\left(\displaystyle 17!17!6!\right)}
{40!}
\displaystyle\prod_{i,j}\frac{1}{x_{ij}!},\ \ \ \Bx \in {\cal F}_{A\Bx^o}.
\]
Then we have the exact conditional $p$-value 
\[
 p = \sum_{\Bx \in {\cal F}_{A\Bx^o}}g(\Bx)h(\Bx) =  0.07035480,
\]
where the test function is 
\[
 g(\Bx) = \left\{\begin{array}{ll}
1, & \chi^2(\Bx) \geq 8.6687,\\
0, & \mbox{otherwise}.
\end{array}
\right.
\]
As a result, we cannot reject H$_0$ at
      significance level $0.05$, which is the same result to strategy
      (a). 

\begin{Example}
\label{example:2x3}
The following toy example should help the reader in understanding the
 method. Let
consider the $2\times 3$ contingency table with the row sums and the
column sums given as follows.
\[
\begin{array}{|c|c|c|c}\cline{1-3}
x_{11} & x_{12} & x_{13} & 3\\ \cline{1-3}
x_{21} & x_{22} & x_{23} & 2\\ \cline{1-3}
\multicolumn{1}{c}{2} & \multicolumn{1}{c}{2} &\multicolumn{1}{c}{1} &
\multicolumn{1}{c}{5} 
\end{array}
\]
There are $5$ elements in the fiber as
\[\begin{array}{rcl}
 {\cal F}_{(3,2,2,2,1)} & = & \left\{\ 
\begin{array}{|c|c|c|}\hline
2 & 1 & 0\\ \hline
0 & 1 & 1\\ \hline
\end{array}\ ,
\hspace*{2mm}
\begin{array}{|c|c|c|}\hline
2 & 0 & 1\\ \hline
0 & 2 & 0\\ \hline
\end{array}\ ,
\hspace*{2mm}
\begin{array}{|c|c|c|}\hline
1 & 2 & 0\\ \hline
1 & 0 & 1\\ \hline
\end{array}\ ,
\hspace*{2mm}
\begin{array}{|c|c|c|}\hline
1 & 1 & 1\\ \hline
1 & 1 & 0\\ \hline
\end{array}\ ,
\hspace*{2mm}
\begin{array}{|c|c|c|}\hline
0 & 2 & 1\\ \hline
2 & 0 & 0\\ \hline
\end{array}\ 
\right\}\\
& = & \{\Bx_1,\ \Bx_2,\ \Bx_3,\ \Bx_4,\ \Bx_5\}.
\end{array}
\]
The fitted value under the ${\cal M}_{indp}$ is 
$
\begin{array}{|c|c|c|}\hline
1.2 & 1.2 & 0.6\\ \hline
0.8 & 0.8 & 0.4\\ \hline
\end{array}$. Then 
the Pearson goodness-of-fit $\chi^2$ for each element is calculated as
\[
 (\chi^2(\Bx_1),\ \chi^2(\Bx_2),\ \chi^2(\Bx_3),\ \chi^2(\Bx_4),\ 
 \chi^2(\Bx_5))
= (2.917,\ 5,\ 2.917,\ 0.833,\ 5).
\]
The conditional probabilities
\[
 h(\Bx) = \frac{3!2!2!2!}{5!}\prod_{i,j}\frac{1}{x_{ij}!} = 
\frac{2}{5}\prod_{i,j}\frac{1}{x_{ij}!}
\]
for each element are calculated
as
\[
 (h(\Bx_1),\ h(\Bx_2),\ h(\Bx_3),\ h(\Bx_4),\ h(\Bx_5)) = (0.2,\ 0.1,\
 0.2,\ 0.4,\ 0.1).
\]
Therefore the conditional $p$-value for $\Bx_4$ is $1.0$, that for
$\Bx_1$ or $\Bx_3$ is $0.6$, and that for $\Bx_2$ or $\Bx_5$ is $0.2$.
\end{Example}

\begin{Remark}
We briefly mention the generalization of the above method to general problems
 and models. First important point is the existence of the
 minimal sufficient statistics in the form of (\ref{eqn:t=Ax}). 
It is known that, 
for the {\it exponential family}, well-known family of the distribution, 
minimal sufficient statistics exist, and 
for a special case of the exponential family, called the {\it toric
 model}, minimal sufficient statistics of the form (\ref{eqn:t=Ax})
 exist. 
The toric model is relatively new concept arising in the field of the
 computational algebraic statistics and is 
defined from the configuration matrix $A = (a_{ij}) \in
 \mathbb{Z}_{\geq 0}^{d\times \nu}$ as follows. 
For the $j$th column vector $\Ba_j = (a_{1j},\ldots,a_{dj})$ of $A$, $j
 \in [\nu]$, define the monomial
\[
 \Btheta^{\Ba_j} = \prod_{i = 1}^d\theta_i^{a_{ij}},\ \ j \in [\nu].
\]
Then the toric model of $A$ is the image of the orthant $\mathbb{R}_{>
 0}^d$ under the map
\[
 f:\ \mathbb{R}^d\ \rightarrow\ \mathbb{R}^{\nu},\ \ \Btheta \mapsto
 \frac{1}{\sum_{j=1}^{\nu}\Btheta^{\Ba_j}}(\Btheta^{\Ba_1},\ldots,\Btheta^{\Ba_{\nu}}).
\]
See Chapter 1.2 of \cite{Pacher-Sturmfels} for detail. 
The toric model specified by the configuration matrix $A \in
 \mathbb{Z}_{\geq 0}^{d\times \nu}$ is also written by
\[
 {\cal M}_A = \{\Btheta = (\theta_i) \in \Delta_{\nu-1}:\ \log \Btheta \in
 \mbox{rowspan}(A)\},
\]
where $\mbox{rowspan}(A) = \mbox{image}(A')$ 
is the linear space spanned by the rows of
 $A$, and $\log\Btheta = (\log \theta_1,\ \ldots,\ \log\theta_{\nu})'$,
 where ${}'$ is a transpose. In
 statistical fields, this is called a {\it log-linear model}. In fact,
 for example of the independence model ${\cal M}_{indp}$ of $2\times 3$
 tables, that is a log-linear model, the parametric description
 $\theta_{ij} = \alpha_i\beta_j$ can
 be written as
\[
\left(\begin{array}{c}
\log\theta_{11}\\
\log\theta_{12}\\
\log\theta_{13}\\
\log\theta_{21}\\
\log\theta_{22}\\
\log\theta_{23}\\
\end{array}
\right) = \left(\begin{array}{ccccc}
1 & 0 & 1 & 0 & 0\\
1 & 0 & 0 & 1 & 0\\
1 & 0 & 0 & 0 & 1\\
0 & 1 & 1 & 0 & 0\\
0 & 1 & 0 & 1 & 0\\
0 & 1 & 0 & 0 & 1
\end{array}
\right)\left(\begin{array}{c}
\alpha_1\\
\alpha_2\\
\beta_1\\
\beta_2\\
\beta_3\\
\end{array}
\right)\ .
\]
The conditional probability function, i.e., 
the generalization of the hypergeometric distribution $h(\Bx)$ in
 (\ref{eqn:hyper-geo})
is as follows. 
For the model specified by the configuration matrix $A$, the conditional
 probability function for given sufficient statistic $A\Bx^o$ is
\[
 P(\BX = \Bx\ |\ A\BX = A\Bx^o) = 
C_{A\Bx^o}^{-1}
\frac{1}{\displaystyle\prod_{i \in
 [\nu]}x_i!},
\]
where
\begin{equation}
 C_{A\Bx^o} = 
\displaystyle\sum_{\By \in {\cal F}_{A\Bx^o}}
\frac{1}{\displaystyle\prod_{i \in [\nu]}y_i!}
\label{eqn:constant-general}
\end{equation}
is a normalizing constant.
Based on this conditional probability function, we can calculate the
 conditional $p$-values by (\ref{eqn:exact-p-value}).

Finally, we note an optimality of the method briefly. 
The conditional procedure mentioned above 
is justified if we consider the
 hypothetical testing to the class
 of {\it similar tests} and the minimal sufficient statistics is {\it
 complete}. For the class of the exponential family, it is known that
 the minimal sufficient statistic is {\it complete}. See Chapter 4.3 of
 \cite{Lehmann-TSH} for detail. 
\end{Remark}

\paragraph*{(c) Estimate the $p$-value by the Monte Carlo method.}
The two strategies to evaluate $p$-values we have considered, asymptotic
evaluation and exact computation, have both advantages and
disadvantages. The asymptotic evaluations relying on the asymptotic
$\chi^2$ distribution are easy to carry out, especially by various
packages in softwares such as \verb|R|. However, poor
fitting to the asymptotic distribution can not be ignorable for
sparse or unbalanced data even with relatively large sample sizes.
The exact calculation of the conditional $p$-values is the best method
if it is possible to carry out. 
In fact, various exact methods and algorithms are considered for
problems of various types of the contingency tables, statistical
models and test statistics. See the survey paper \cite{agresti-1992} for this
field. However, for large size samples, 
the cardinality of the fiber $|{\cal F}_{A\Bx^o}|$ can exceed billions,
making exact computations difficult to be carried out. 
In fact, it is known that the cardinality of a fiber increases
exponentially in the sample size $n$. (An approximation for the
cardinality of a fiber is given by \cite{Gail-Mantel}.)  
For these cases, the Monte Carlo methods can be effective.

The Monte Carlo methods estimate the $p$-values as follows. To compute
the conditional $p$-value (\ref{eqn:exact-p-value}), generate samples
$\Bx_1,\ldots,\Bx_N$ from the null distribution
$h(\Bx)$. Then the $p$-value is estimated as $\hat{p} = \sum_{i = 1}^N
g(\Bx_i) / N$, that is an unbiased estimate of the $p$-value. We can set
$N$ according to the performance of our computer. As an advantage of the
Monte Carlo method, we can also estimate the accuracy, i.e., {\it
variance} of the estimate. For example, a conventional $95\%$
confidence interval of $p$, $\hat{p} \pm
1.96\sqrt{\hat{p}(1-\hat{p})/N}$, is frequently used.  
The problem here is how to generate samples from the null
distribution. We consider {\it Markov chain Monte Carlo methods}, often
abbreviated as the MCMC methods, in this paper. 

Following MCMC methods setup, we construct an ergodic Markov chain 
on the fiber ${\cal F} = {\cal F}_{A\Bx^o}$ whose 
stationary distribution is prescribed, given by (\ref{eqn:hyper-geo}). 
Let
the elements of ${\cal F}$ be numbered as 
\[
 {\cal F} =
\{\Bx_1,\ldots,\Bx_s\}.
\]
We write the null distribution on ${\cal F}$ as
\[
 \Bpi = (\pi_1,\ldots,\pi_s) = (h(\Bx_1),\ldots,h(\Bx_s)).
\]
Here, by standard notation, we treat
$\Bpi$ as a row vector. We write the transition probability matrix of
the Markov chain $\{Z_t, t \in \mathbb{Z}_{\geq 0}\}$ over ${\cal F}$ as
$Q = (q_{ij})$, i.e., we define
\[
 q_{ij} = P(Z_{t+1} = \Bx_j\ |\ Z_t = \Bx_i).
\]
Then a probability distribution $\Btheta \in \Delta_{s-1}$
is called a stationary distribution if it satisfies $\Btheta =
\Btheta Q$. The stationary distribution uniquely exists if the Markov chain
is irreducible, (i.e., connected in this case) and aperiodic. Therefore
for the connected and aperiodic Markov chain, starting from an arbitrary
state $Z_0 = \Bx_i$, the distribution of $Z_t$ for large $t$ is close to
its stationary distribution. If we can construct a connected and aperiodic
Markov chain with the stationary distribution $\Bpi$, by
running the Markov chain and discarding a large number $t$ of initial
steps (called {\it burn-in steps}), we can treat $Z_{t+1}, Z_{t+2},
\ldots$ to be samples from the null distribution $\Bpi$ and use
them to estimate $p$-values. Then the problem becomes {\it how to
construct a connected and aperiodic Markov chain with the stationary
distribution as the null distribution $\Bpi$ over ${\cal F}$}. Among
these conditions, the conditions for the stationary distribution can be
solved easily. Once we construct an arbitrary connected chain over
${\cal F}$, we can modify its stationary distribution to the given null
distribution $\Bpi$ as follows.
\begin{Theorem}[Metropolis-Hastings algorithm]
\label{thm:M-H-algo}
Let $\Bpi$ be a probability distribution on ${\cal F}$. Let $R =
 (r_{ij})$ be the transition probability matrix of a connected,
 aperiodic and symmetric Markov chain over ${\cal F}$. Then the
transition probability matrix $Q = (q_{ij})$ defined by
\[
\begin{array}{rcl}
 q_{ij} & = & r_{ij}\min\displaystyle\left(1,
				      \frac{\pi_j}{\pi_i}\right),\ i\neq
 j\\
q_{ii} & = & 1 - \displaystyle\sum_{j \neq i}q_{ij}
\end{array}
\]
satisfies $\Bpi = \Bpi Q$.
\end{Theorem}
This theorem is a special case of \cite{Hastings-1970}. Though the
symmetry assumption ($r_{ij} = r_{ji}$) can be removed easily, we
only consider symmetric $R$ for simplicity. The proof of this theorem is
easy and is omitted. See \cite{Hastings-1970} or Chapter 4.1 of
\cite{dojo-en}, for example. Instead, we consider the algorithm for
data of small size.
\begin{Example}
\label{example:2x3-cont}
Consider the small example in Example \ref{example:2x3}. As we have
 seen, the fiber is
\[
 {\cal F} = \{\Bx_1,\ \Bx_2,\ \Bx_3,\ \Bx_4,\ \Bx_5\}
\]
and the null distribution is 
\[
 \Bpi = (\pi_1,\ldots,\pi_5) = (h(\Bx_1),\ldots,h(\Bx_5)) = (0.2, 0.1,
 0.2, 0.4, 0.1).
\]
Using the Markov basis we consider in the next section, we can construct
 a connected, aperiodic and symmetric Markov chain with the 
transition
 probability matrix
\begin{equation}
 R = \left(\begin{array}{ccccc}
1/2 & 1/6 & 1/6 & 1/6 & 0\\
1/6 & 2/3 & 0 & 1/6 & 0\\
1/6 & 0 & 1/2 & 1/6 & 1/6\\
1/6 & 1/6 & 1/6 & 1/3 & 1/6\\
0 & 0 & 1/6 & 1/6 & 2/3
\end{array}
\right)\ .
\label{eqn:R-in-example} 
\end{equation}
Following Theorem \ref{thm:M-H-algo}, 
we modify the Markov chain to
 have the transition
 probability matrix
\[
Q = \left(\begin{array}{ccccc}
7/12 & 1/12 & 1/6 & 1/6 & 0\\
1/6 & 2/3 & 0 & 1/6 & 0\\
1/6 & 0 & 7/12 & 1/6 & 1/12\\
1/12 & 1/24 & 1/12 & 3/4 & 1/24\\
0 & 0 & 1/6 & 1/6 & 2/3
\end{array}
\right)\ .
\]
We can check that the eigenvector from the left of $Q$ with the
 eigenvalue $1$ is $\Bpi$. We can also check that each row vector of
 $Q^T$ for large $T$ converges to $\Bpi$.

\end{Example}
An important advantage of the Markov chain Monte Carlo method is that it
does not require the explicit evaluation of the normalizing constant of
the null distribution.
As is shown in Theorem \ref{thm:M-H-algo}, we only need to know
$\Bpi$ up to a multiplicative constant, because the normalizing
constant, (\ref{eqn:constant-general}) in the general form, canceled in
the ratio $\pi_j/\pi_i$. With Theorem \ref{thm:M-H-algo}, the remaining
problem is to construct an arbitrary connected and aperiodic Markov
chain over ${\cal F}$, that is solved by the Gr\"obner basis theory.

\section{Markov bases and ideals}
As stated in the previous section, 
the main task for estimating $p$-values thanks to MCMC methods is to 
construct a connected and aperiodic Markov
chain over ${\cal F} = {\cal F}_{A\Bx^o}$ with stationary distribution
given by (\ref{eqn:hyper-geo}). 
Here, $A \in \mathbb{Z}^{d\times \nu}$ is
a given configuration matrix, $\Bx^o \in \mathbb{Z}_{\geq 0}^{\nu}$ is
the observed contingency table and ${\cal F}_{A\Bx^o}$, a
$A\Bx^o$-fiber,
 is the set of all
contingency tables with the same value of the minimal sufficient
statistics to $\Bx^o$,
\[
{\cal F}_{A\Bx^o} = \{\Bx \in \mathbb{Z}_{\geq 0}^{\nu}\ :\ A\Bx = A\Bx^o\}. 
\]
We write the integer kernel of $A$ as
\[
 {\rm Ker}_{\mathbb{Z}}(A) = {\rm Ker}(A) \cap \mathbb{Z}^{\nu}
= \{\Bz \in \mathbb{Z}^{\nu}\ :\ A\Bz = \Bzero\}.
\]
An element of $ {\rm Ker}_{\mathbb{Z}}(A)$ is called a {\it move}. 
Note that $\Bx - \By \in {\rm Ker}_{\mathbb{Z}}(A)$ if and only if $\Bx,
\By \in {\cal F}_{A\Bx}$.
Then for a given subset 
${\cal B} \subset {\rm Ker}_{\mathbb{Z}}(A)$ and 
$\Bt \in \mathbb{Z}_{\geq 0}^d$, 
we can define
      undirected graph $G_{\Bt, {\cal B}} = (V, E)$ by
\[
 V = {\cal F}_{\Bt},\ \ \ 
E = \{(\Bx, \By)\ :\ \Bx - \By \in {\cal B}\ \mbox{or}\  
\By - \Bx \in {\cal B}
\}.
\]
\begin{Definition}[A Markov basis]
${\cal B} \subset {\rm Ker}_{\mathbb{Z}}(A)$ is a Markov basis for $A$
      if $G_{\Bt, {\cal B}}$ is connected for arbitrary $\Bt \in
      \mathbb{Z}_{\geq 0}^{d}$. 
\end{Definition}
Once we obtain a Markov basis ${\cal B}$ for $A$, we can construct a
      connected Markov chain over ${\cal F}_{A\Bx^o}$ easily as
      follows. 
For each state $\Bx \in {\cal F}_{A\Bx^o}$,
randomly choose a move $\Bz \in {\cal B}$ and a sign $\varepsilon \in
	     \{-1,1\}$ 
and consider $\Bx + \varepsilon\Bz$. 
If $\Bx + \varepsilon\Bz \in {\cal F}_{A\Bx^o}$, then $\Bx
	     + \varepsilon\Bz$ is
	     the next state, otherwise stay at $\Bx$. Then we have the
	     connected Markov chain over $A\Bx^o$. 
We see these arguments in an example.
\begin{Example}
\label{example:2x3-cont2}
Again we consider a small data of Example \ref{example:2x3}, where the
 fiber is redisplayed below.
\[\begin{array}{rcl}
 {\cal F}_{(3,2,2,2,1)} & = & \left\{\ 
\begin{array}{|c|c|c|}\hline
2 & 1 & 0\\ \hline
0 & 1 & 1\\ \hline
\end{array}\ ,
\hspace*{2mm}
\begin{array}{|c|c|c|}\hline
2 & 0 & 1\\ \hline
0 & 2 & 0\\ \hline
\end{array}\ ,
\hspace*{2mm}
\begin{array}{|c|c|c|}\hline
1 & 2 & 0\\ \hline
1 & 0 & 1\\ \hline
\end{array}\ ,
\hspace*{2mm}
\begin{array}{|c|c|c|}\hline
1 & 1 & 1\\ \hline
1 & 1 & 0\\ \hline
\end{array}\ ,
\hspace*{2mm}
\begin{array}{|c|c|c|}\hline
0 & 2 & 1\\ \hline
2 & 0 & 0\\ \hline
\end{array}\ 
\right\}\\
& = & \{\Bx_1,\ \Bx_2,\ \Bx_3,\ \Bx_4,\ \Bx_5\}.
\end{array}
\]
The integer kernel for the configuration matrix
\begin{equation}
 A = \left(\begin{array}{cccccc}
1 & 1 & 1 & 0 & 0 & 0\\
0 & 0 & 0 & 1 & 1 & 1\\
1 & 0 & 0 & 1 & 0 & 0\\
0 & 1 & 0 & 0 & 1 & 0\\
0 & 0 & 1 & 0 & 0 & 1
\end{array}
\right)
\label{eqn:A-2x3-indp}
\end{equation}
includes moves such as
\[
\Bz_1 = \begin{array}{|r|r|r|}\hline
1 & -1 & 0\\ \hline
-1 & 1 & 0\\ \hline
\end{array}\ ,\ \  
\Bz_2 = \begin{array}{|r|r|r|}\hline
1 & 0 & -1\\ \hline
-1 & 0 & 1\\ \hline
\end{array}\ ,\ \  
\Bz_3 = \begin{array}{|r|r|r|}\hline
0 & 1 & -1 \\ \hline
0 & -1 & 1 \\ \hline
\end{array}\ ,\ \  
\begin{array}{|r|r|r|}\hline
2 & -1 & -1\\ \hline
-2 & 1 & 1\\ \hline
\end{array}\ ,\ \ldots. 
\]
From these, we consider some sets of moves. 
If we consider 
${\cal B}_1 = \{\Bz_1\}$, corresponding undirected graph
 $G_{(3,2,2,2,1), {\cal B}_1}$ is given in Figure \ref{fig:graph-B-2x3}(a),
 which is not connected. Therefore ${\cal B}_1$ is not a Markov
 basis. If we consider 
${\cal B}_2 = \{\Bz_1, \Bz_2\}$, corresponding undirected graph
 $G_{(3,2,2,2,1), {\cal B}_2}$ is given in Figure
 \ref{fig:graph-B-2x3}(b), which is connected. However, ${\cal B}_2$ is
 also not a Markov basis, because there exists $\Bt \in \mathbb{Z}_{\geq
 0}^5$ where $G_{\Bt, {\cal B}_2}$ is not connected. An example of such
 $\Bt$ is $\Bt = (1,1,0,1,1)$, with the corresponding $\Bt$-fiber is a
 two-element set
\begin{equation}
 {\cal F}_{(1,1,0,1,1)} = \left\{
\ \begin{array}{|r|r|r|}\hline
0 & 1 & 0\\ \hline
0 & 0 & 1\\ \hline
\end{array}\ ,\ \ 
\begin{array}{|r|r|r|}\hline
0 & 0 & 1\\ \hline
0 & 1 & 0\\ \hline
\end{array}\ 
\right\}\ .
\label{eqn:F-11011}
\end{equation}
The above example shows that a Markov basis includes $\Bz_3$ to connect
 the two elements above. In fact, ${\cal B} = \{\Bz_1, \Bz_2, \Bz_3\}$
 is a Markov basis for this $A$, with the corresponding 
undirected graph 
 $G_{(3,2,2,2,1), {\cal B}_3}$ in Figure
 \ref{fig:graph-B-2x3}(c).
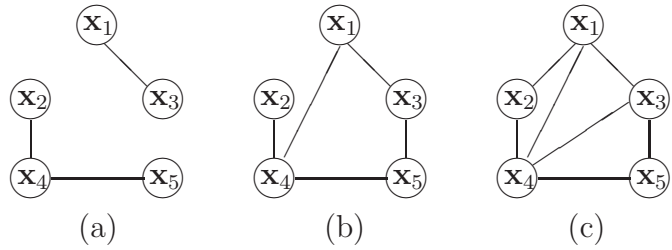
\begin{figure}[htbp]
\begin{center}
\begin{picture}(70,80)(0,0)
\put(35,63){$\Bx_1$}
\put(40,65){\circle{15}}
\put(10,35){$\Bx_2$}
\put(15,37){\circle{15}}
\put(60,35){$\Bx_3$}
\put(65,37){\circle{15}}
\put(10,5){$\Bx_4$}
\put(15,7){\circle{15}}
\put(60,5){$\Bx_5$}
\put(65,7){\circle{15}}
\put(59,42){\line(-1,1){16}}
\put(15,15){\line(0,1){14}}
\put(23,7){\line(1,0){35}}
\put(33,-15){(a)}
\end{picture}
\hspace*{5mm}
\begin{picture}(70,80)(0,0)
\put(35,63){$\Bx_1$}
\put(40,65){\circle{15}}
\put(10,35){$\Bx_2$}
\put(15,37){\circle{15}}
\put(60,35){$\Bx_3$}
\put(65,37){\circle{15}}
\put(10,5){$\Bx_4$}
\put(15,7){\circle{15}}
\put(60,5){$\Bx_5$}
\put(65,7){\circle{15}}
\put(59,42){\line(-1,1){16}}
\put(19,15){\line(1,2){21}}
\put(15,15){\line(0,1){14}}
\put(65,15){\line(0,1){14}}
\put(23,7){\line(1,0){35}}
\put(34,-15){(b)}
\end{picture}
\hspace*{5mm}
\begin{picture}(70,80)(0,0)
\put(35,63){$\Bx_1$}
\put(40,65){\circle{15}}
\put(10,35){$\Bx_2$}
\put(15,37){\circle{15}}
\put(60,35){$\Bx_3$}
\put(65,37){\circle{15}}
\put(10,5){$\Bx_4$}
\put(15,7){\circle{15}}
\put(60,5){$\Bx_5$}
\put(65,7){\circle{15}}
\put(21,42){\line(1,1){16}}
\put(59,42){\line(-1,1){16}}
\put(19,15){\line(1,2){21}}
\put(21,12){\line(3,2){36}}
\put(15,15){\line(0,1){14}}
\put(65,15){\line(0,1){14}}
\put(23,7){\line(1,0){35}}
\put(34,-15){(c)}
\end{picture}
\end{center}
\caption{Undirected graphs for ${\cal B}_1, {\cal B}_2, {\cal B}_3$ for
 $\Bt = (3,2,2,2,1)$.}
\label{fig:graph-B-2x3},
\end{figure}
The transition probability matrix (\ref{eqn:R-in-example}) in
 Example \ref{example:2x3-cont} corresponds to a Markov chain
 constructed from ${\cal B}_3$ as ``in each step, choose $3$ elements in
 ${\cal B}_3$
 and its sign $\{-1,1\}$ with equal probabilities''.
\end{Example}
At first sight, we may feel the cases such as (\ref{eqn:F-11011}) are
 trivial and may imagine that ``if we only consider the cases with $\Bt \in
 \mathbb{Z}_{>0}^d$, i.e., cases with strictly positive minimal sufficient
 statistics (that may be realistic situations in the actual data analysis),
 it is easy to connect the fiber ${\cal F}_{\Bt}$''. However, it is not
 so. We will see an example where complicated moves are needed even for
 the fiber with positive $\Bt$. 

The connection between the Markov basis and a {\it toric ideal} of a
polynomial ring by \cite{Diaconis-Sturmfels-1998} is as follows. Let 
$k[\Bu] = k[u_{1}, u_{2},\ldots,u_{\nu}]$ denote the ring of polynomials
in $\nu$ variables over a field $k$. Let 
a contingency table $\Bx \in \mathbb{Z}_{\geq
0}^{\nu}$ be mapped to the
      monomial $\Bu^{\Bx} \in k[\Bu]$, and a move, i.e., 
an element of the integer kernel $\Bz =
      \Bz^+ - \Bz^- \in {\rm
      Ker}_{\mathbb{Z}}(A)$, be mapped to the binomial $\Bu^{\Bz^+} -
      \Bu^{\Bz^-} \in k[\Bu]$. For the case of the independence model
      for the $3\times 3$ contingency
      tables, examples of these correspondences are as follows.
\[
\begin{array}{|r|r|r|cl}\cline{1-3}
11 & 5 & 2 & & \\ \cline{1-3}
4 & 9 & 1 & \Longleftrightarrow &
u_{11}^{11}u_{12}^5u_{13}^2u_{21}^4u_{22}^9u_{23}u_{31}^2u_{32}^3u_{33}^3
\\ \cline{1-3} 
2 & 3 & 3 & & \\ \cline{1-3}
\multicolumn{5}{c}{}\\ \cline{1-3}
2 & -1 & -1 & & \\ \cline{1-3}
-3 & 1 & 2 & \Longleftrightarrow & u_{11}^{2}u_{22}u_{23}^2u_{31} -
      u_{12}u_{13}u_{21}^3u_{33}\\ \cline{1-3} 
1 & 0 & -1 & & \\ \cline{1-3}
 \end{array} 
\]
The binomial ideal in $k[\Bu]$ generated by the set of
binomials corresponding to the set of moves for $A$, 
\[
 I_A = \left<\left\{
\Bu^{\Bz^+} - \Bu^{\Bz^-} :\ \Bz^+ - \Bz^- \in {\rm
      Ker}_{\mathbb{Z}}(A)
\right\}\right>, 
\]
is the {\it the toric ideal} of configuration $A$. 

\begin{Theorem}[Theorem 3.1 of \cite{Diaconis-Sturmfels-1998}]
\label{thm:thm3.1-DS}
${\cal B} = \{\Bz_1,\ldots,\Bz_L\} \subset {\rm Ker}_{\mathbb{Z}}(A)$ is
      a Markov basis for $A$ if and only if $\{\Bu^{\Bz_i^+} -
      \Bu^{\Bz_i^-},\ i = 1,\ldots,L\}$ generates $I_A$.
\end{Theorem}
A proof of Theorem \ref{thm:thm3.1-DS} is given in the original paper 
\cite{Diaconis-Sturmfels-1998}. We can also find more detailed proof in
Chapter 4 of \cite{dojo-en}. In these proofs, the sufficiency and the
necessity are shown by induction on some integer. In the proof of
sufficiency, this integer represents the number of steps of the chain,
and the argument is straightforward. On the other hand, in the proof of
necessity, this integer represents the number of 
terms in the expansion that we want to show in
the proof, and is not necessarily equal to the number of steps of the
chain. Theorem \ref{thm:thm3.1-DS} shows a non-trivial result on this point.

To calculate a Markov basis for a given configuration matrix $A$, 
we can use the {\it elimination theory}. For this purpose, we also
prepare variables $\Bv = \{v_1,\ldots,v_{d}\}$ for the minimal
sufficient statistic $\Bt$ and consider the polynomial ring
$k[\Bv] = k[v_1,\ldots,v_{d}]$. The relation $\Bt = A\Bx$ can be
expressed by the homomorphism
\[\begin{array}{rcl}
\psi_A\ :\ k[\Bu] & \rightarrow & k[\Bv]\\
u_j & \mapsto & v_1^{a_{1j}}v_2^{a_{2j}}\cdots v_{d}^{a_{dj}}.
\end{array}
\]
Then the toric ideal $I_A$ is also expressed as $I_A = {\rm
	     Ker}(\psi_A)$.
We now have the following. 
\begin{Corollary}[Theorem 3.2 of \cite{Diaconis-Sturmfels-1998}]
\label{cor:elim-thm}
Let $I_A^{*}$ be the ideal of $k[\Bu, \Bv]$ given by
\[
 I_A^{*} = \left<-\psi_A(u_j) + u_j,\ j = 1,\ldots,\nu \right> \subset
	     k[\Bu, \Bv].
\]
Then we have $I_A = I_A^{*} \cap k[\Bu]$.
\end{Corollary}
Corollary \ref{cor:elim-thm} suggests that we can obtain a generator of
$I_A$ as its Gr\"obner basis for an appropriate term order called an
elimination order. 
For an ideal $J \in k[\Bu]$ and a term order $\prec$, a set of
polynomials $\{g_1,\ldots,g_s\},\ g_1,\ldots,g_s \in J$, is called a 
 Gr\"obner basis of $I$ with respect to a term order $\prec$, if
 $\{{\rm in}_{\prec}(g_1), \ldots, {\rm in}_{\prec}(g_s)\}$ generates 
an initial ideal of $J$ defined by ${\rm in}_{\prec}(J) = \left<\{{\rm
in}_{\prec}(f)\ :\ 0 \neq f \in J\}\right>$. Here we write ${\rm
in}_{\prec}(f)$ as an initial term of $f$ with respect to a term order
$\prec$. For more theories and results on Gr\"obner bases, see textbooks
such as \cite{Cox-Little-Oshea}. The elimination theory is one of the
useful applications of Gr\"obner bases and is used for our problem as follows. 
For the reduced Gr\"obner basis $G^*$ of
	     $I_A^*$ for any term order satisfying $\{v_1,\ldots,v_d\} \succ
	     \{u_1,\ldots,u_{\nu}\}$, $G^{*} \cap k[\Bu]$ is a reduced
	     Gr\"obner basis of $I_A$. Because the Gr\"obner basis is a
 generator of $I_A$, we can obtain a Markov basis for $A$ as the reduced
 Gr\"obner basis in this way.

The computations of Gr\"obner bases can be carried out by various
algebraic softwares such as Macaulay2 (\cite{Macaulay2}), SINGULAR
(\cite{Singular}),
CoCoA (\cite{CoCoA}), Risa/Asir (\cite{Risa/Asir}) and 4ti2
(\cite{4ti2}). Here, we show some computations by Macaulay2, because we
can also rapidly use it online at the website
\footnote{Macaulay2 online: {\tt http://habanero.math.cornell.edu.3690}}. 
We start with a simple example.
\begin{Example}
\label{example-2x3-cont4}
In Example \ref{example:2x3-cont2}, we give a Markov basis for the
 independence model for $2\times 3$ contingency tables without any proof or
 calculations. Here we check that the set of $3$ moves
\[
\left\{\ 
\Bz_1 = \begin{array}{|r|r|r|}\hline
1 & -1 & 0\\ \hline
-1 & 1 & 0\\ \hline
\end{array}\ ,\ \  
\Bz_2 = \begin{array}{|r|r|r|}\hline
1 & 0 & -1\\ \hline
-1 & 0 & 1\\ \hline
\end{array}\ ,\ \  
\Bz_3 = \begin{array}{|r|r|r|}\hline
0 & 1 & -1 \\ \hline
0 & -1 & 1 \\ \hline
\end{array}\ 
\right\}
\]
constitute a Markov basis for $A$ given in (\ref{eqn:A-2x3-indp}).
In other words, we check that 
 the corresponding toric ideal $I_A$ is generated by $3$ binomials
\begin{equation}
 \{u_{11}u_{22} - u_{12}u_{21},\ u_{11}u_{23}-u_{13}u_{21},\
 u_{12}u_{23}-u_{13}u_{22}\}. 
\label{eqn:2x3-MB}
\end{equation}
Following Corollary \ref{cor:elim-thm}, we
prepare the variable $\Bv = (v_1,\ldots,v_5)$ for the row sums and
 column sums of $\Bx$ as
\[
\begin{array}{|r|r|r|r}\cline{1-3}
x_{11} & x_{12} & x_{13} & v_1\\ \cline{1-3}
x_{21} & x_{22} & x_{23} & v_2\\ \cline{1-3}
\multicolumn{1}{c}{v_3} & \multicolumn{1}{c}{v_4} &
\multicolumn{1}{c}{v_5} & 
\end{array}
\]
and consider the homomorphism
\[
 \begin{array}{lll}
u_{11} \mapsto v_1v_3, & u_{12} \mapsto v_1v_4, & 
u_{13} \mapsto v_1v_5,\\
u_{21} \mapsto v_2v_3, & u_{22} \mapsto v_2v_4, & 
u_{23} \mapsto v_2v_5.
 \end{array}
\]
Then under the elimination order $\Bv \succ \Bu$, compute the reduced
 Gr\"obner basis of the toric ideal
\[
 I_{A}^* = \left< -v_1v_3 + u_{11},\ -v_1v_4 + u_{12}, \ldots, -v_2v_5 + u_{23}
\right>.
\]
These calculations are done by Macaulay2 as follows.
{\footnotesize
\begin{verbatim}
i1 : R=QQ[v1,v2,v3,v4,v5,u11,u12,u13,u21,u22,u23,MonomialOrder=>{5,6}]

o1 = R

o1 : PolynomialRing

i2 : I=ideal(-v1*v3+u11,-v1*v4+u12,-v1*v5+u13,-v2*v3+u21,-v2*v4+u22,-v2*v5+u23)

o2 = ideal (- v1*v3 + u11, - v1*v4 + u12, - v1*v5 + u13, - v2*v3 + u21, - v2*v4
     --------------------------------------------------------------------------
     + u22, - v2*v5 + u23)

o2 : Ideal of R

i3 : G=gb(I); g=gens(G)

o4 = | u13u22-u12u23 u13u21-u11u23 u12u21-u11u22 v4u23-v5u22 v4u13-v5u12
     --------------------------------------------------------------------------
     v3u23-v5u21 v3u22-v4u21 v3u13-v5u11 v3u12-v4u11 v1u23-v2u13 v1u22-v2u12
     --------------------------------------------------------------------------
     v1u21-v2u11 v2v5-u23 v1v5-u13 v2v4-u22 v1v4-u12 v2v3-u21 v1v3-u11 |

             1       18
o4 : Matrix R  <--- R

i5 : selectInSubring(1,g)

o5 = | u13u22-u12u23 u13u21-u11u23 u12u21-u11u22 |

             1       3
o5 : Matrix R  <--- R
\end{verbatim}
}
The output \verb|o4| shows the reduced Gr\"obner basis of $I_A^*$ under the
 elimination (reverse lexicographic) order $\Bv \succ \Bu$, and the
 output \verb|o5| shows the reduced Gr\"obner basis of $I_A$, which we
 can use as a Markov basis. We have now checked a Markov basis
 (\ref{eqn:2x3-MB}). 
\end{Example}

From the Markov basis (\ref{eqn:2x3-MB}), we may imagine that the set of
moves corresponding to the binomials
\[
\{ u_{ij}u_{i'j'} - u_{ij'}u_{i'j}\ ,\ \ 1 \leq i < i' \leq I,\ 1 \leq j <
j' \leq J\}
\]
forms a Markov basis for the independence model of the $I\times J$
contingency tables, which is actually true. This fact is given and
proved as
Theorem 2.1 of \cite{Aoki-Hara-Takemura-2012}, for example.

\begin{figure}[htbp]
\begin{center}
\includegraphics[width=80mm]{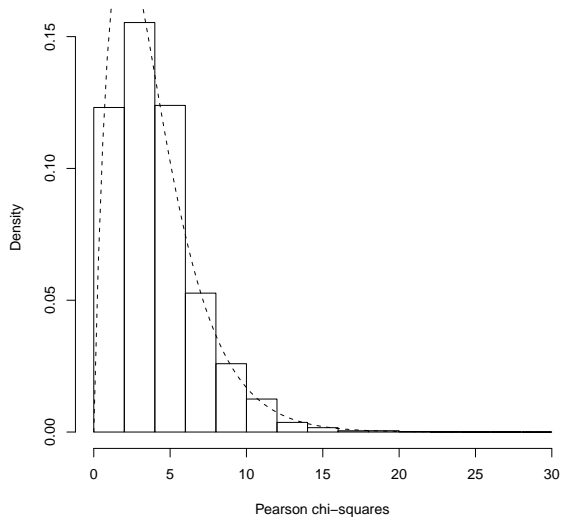}
\caption{A histogram of sampled Pearson $\chi^2$ goodness-of-fit 
for Table \ref{tbl:3x3-example} 
generated by a Markov chain Monte Carlo method. 
The dotted curve is
 the corresponding asymptotic $\chi^2_4$ distribution.}
\label{fig:mcmc-conditional-distribution-Table1}
\end{center}
\end{figure}

Now we are ready to estimate $p$-value for 
our original problem of $3\times 3$ contingency table in Table
\ref{tbl:3x3-example}. The Markov basis for this problem is formed by 
$9$ moves of the above type. Using this Markov basis, 
we calculate the conditional $p$-values for
Table \ref{tbl:3x3-example} by the Markov chain Monte Carlo method.  
For each step of the chain, we choose an element of the Markov basis
randomly, and modify the transition probability by Theorem
\ref{thm:M-H-algo}. We start the chain at the observed table
$\Bx^o$ of Table \ref{tbl:3x3-example}, discard initial $50000$
steps as the burn-in steps, and have $100000$ samples of the Pearson
goodness-of-fit $\chi^2$. 
Figure \ref{fig:mcmc-conditional-distribution-Table1} is a 
histogram of the sampled Pearson
goodness-of-fit $\chi^2$ with the asymptotic $\chi^2_4$ distribution.
In these $100000$ samples, $6681$
samples are larger than or equal to the observed value $\chi^2(\Bx^o) =
8.6687$, then we have the estimate $\hat{p} = 0.06681$. Therefore we
cannot reject $\mbox{H}_0$ at significance level $0.05$, which is the
same result to the other strategies (a) and (b). Though the difference
from the exact value $p = 0.07035480$ from the simulated value 
is slightly larger than the asymptotic estimate
($\hat{p} = 0.0699$), we may increase the accuracy of the estimates by
increasing the sample sizes. 
To compare the three strategies for Table \ref{tbl:3x3-example}, we
compute the upper percentiles of $90\%, 95\%, 99\%, 99.9\%$ for (a)
asymptotic $\chi^2_4$ distribution, (b) exact conditional distribution
and (c) Monte Carlo simulated distribution in Table
\ref{tbl:percentiles-example}.
\begin{table}[htbp]
\begin{center}
\caption{The upper percentiles for three strategies of Pearson
 goodness-of-fit $\chi^2$ for Table \ref{tbl:3x3-example}.}
\label{tbl:percentiles-example}
\begin{tabular}{lllll}\hline
                                       & $90\%$ & $95\%$ & $99\%$ &
 $99.9\%$\\ \hline
(a) Asymptotic $\chi^2_4$ distribution & $7.779$ & $9.488$ & $13.28$ &
		 $18.47$\\
(b) Exact null distribution            & $7.766$ & $9.353$ & $12.78$ &
		 $17.99$\\
(c) Monte Carlo simulated distribution & $7.684$ & $9.287$ & $12.73$ &
		 $18.58$\\ \hline
\end{tabular}
\end{center}
\end{table}

Finally, we give an example for which the structure of the 
Markov basis is
complicated. 
The model we consider is a {\it no three-factor interaction model} for
three-way contingency tables. 
The parametric description of the no three-factor interaction model is
given by
\[
 {\cal M}_{n3} = \{\Btheta \in \Delta\ :\ \theta_{ijk} =
 \alpha_{ij}\beta_{ik}\gamma_{jk} \ \mbox{for some}\ (\alpha_{ij}),
 (\beta_{ik}), (\gamma_{jk})\}.
\]
This is one of the most important statistical models in the statistical
data analysis of three-way contingency tables. 
The minimal sufficient statistics for ${\cal M}_{n3}$ is the
two-dimensional marginals
\[
 \{x_{ij+}\},\ \{x_{i+k}\},\ \{x_{+jk}\},
\]
where we define
\[
 x_{ij+} = \sum_{k =1}^Kx_{ijk},\ \ 
 x_{i+k} = \sum_{j =1}^Jx_{ijk},\ \ 
 x_{+jk} = \sum_{i =1}^Ix_{ijk}.
\]
We only consider $3\times 3\times 3$ case (i.e., $I = J = K = 3$) here. 
Then the configuration matrix $A$ is $27\times 27$ matrix 
written as follows.
\[{\scriptsize
 A = \left(\begin{array}{ccccccccccccccccccccccccccc}
1 & 1 & 1 & 0 & 0 & 0 & 0 & 0 & 0 & 0 & 0 & 0 & 0 & 0 & 0 & 0 & 0 & 0 &
 0 & 0 & 0 & 0 & 0 & 0 & 0 & 0 & 0\\
0 & 0 & 0 & 1 & 1 & 1 & 0 & 0 & 0 & 0 & 0 & 0 & 0 & 0 & 0 & 0 & 0 & 0 &
 0 & 0 & 0 &
 0 & 0 & 0 & 0 & 0 & 0\\
\multicolumn{27}{c}{\vdots}\\
0 & 0 & 0 & 0 & 0 & 0 & 0 & 0 & 0 & 0 & 0 & 0 & 0 & 0 & 0 &
 0 & 0 & 0 &
 0 & 0 & 0 & 0 & 0 & 0 & 1 & 1 & 1\\
1 & 0 & 0 & 1 & 0 & 0 & 1 & 0 & 0 & 0 & 0 & 0 & 0 & 0 & 0 & 0 & 0 & 0 &
 0 & 0 & 0 &
 0 & 0 & 0 & 0 & 0 & 0\\
0 & 1 & 0 & 0 & 1 & 0 & 0 & 1 & 0 & 0 & 0 & 0 & 0 & 0 & 0 & 0 & 0 & 0 & 0 &
 0 & 0 & 0 &
 0 & 0 & 0 & 0 & 0\\
\multicolumn{27}{c}{\vdots}\\
0 & 0 & 0 & 0 & 0 & 0 & 0 & 0 & 0 & 0 & 0 & 0 & 0 & 0 & 0 & 0 & 0 & 0 & 0 &
 0 & 1 & 0 &
 0 & 1 & 0 & 0 & 1\\
1 & 0 & 0 & 0 & 0 & 0 & 0 & 0 & 0 & 
1 & 0 & 0 & 0 & 0 & 0 & 0 & 0 & 0 & 
1 & 0 & 0 & 0 & 0 & 0 & 0 & 0 & 0\\
0 & 1 & 0 & 0 & 0 & 0 & 0 & 0 & 0 & 0 & 
1 & 0 & 0 & 0 & 0 & 0 & 0 & 0 & 0 & 
1 & 0 & 0 & 0 & 0 & 0 & 0 & 0\\
\multicolumn{27}{c}{\vdots}\\
0 & 0 & 0 & 0 & 0 & 0 & 0 & 0 & 1 & 
0 & 0 & 0 & 0 & 0 & 0 & 0 & 0 & 1 & 
0 & 0 & 0 & 0 & 0 & 0 & 0 & 0 & 1
\end{array}
\right)
}\]
For this model we see that the ``simplest moves'', i.e., the moves with
the minimum degree, correspond to the binomials of degree $4$ such as
\begin{equation}
 u_{111}u_{122}u_{212}u_{221} - u_{112}u_{121}u_{211}u_{222},
\label{eqn:basic-move-3way}
\end{equation}
which is called a {\it basic move}. 
There are $9$ such moves for the case of $3\times 3\times 3$ tables. 
Unfortunately, however, the set of these $9$ moves does not become a
Markov basis. To see this consider the following example.

\begin{Example}
Consider the $3\times 3\times 3$ contingency tables with the fixed
 two-dimensional marginals 
\begin{equation}
 (x_{ij+}) = (x_{i+k}) = (x_{+jk}) = (2,1,1,1,2,1,1,1,2)'.
\label{eqn:3x3x3-marginal}
\end{equation}
We write a $3\times 3\times 3$ table as follows.
\[
\begin{array}{|c|}\hline
x_{111}\ x_{112}\ x_{113}\\ x_{121}\ x_{122}\ x_{123}\\ x_{131}\
 x_{132}\ x_{133}\\ \hline
\end{array}
\hspace*{1mm}
\begin{array}{|c|}\hline
x_{211}\ x_{212}\ x_{213}\\ x_{221}\ x_{222}\ x_{223}\\ x_{231}\
 x_{232}\ x_{233}\\ \hline
\end{array}
\hspace*{1mm}
\begin{array}{|c|}\hline
x_{311}\ x_{312}\ x_{313}\\ x_{321}\ x_{322}\ x_{323}\\ x_{331}\
 x_{332}\ x_{333}\\ \hline
\end{array}
\]
Then the fixed marginals (\ref{eqn:3x3x3-marginal}) are displayed as
\[
\begin{array}{|@{}c@{}c@{}c@{}|c}\cline{1-3}
 x_{111} & x_{112}& x_{113}& 2\\
 x_{121} & x_{122}& x_{123}& 1\\
 x_{131} & x_{132}& x_{133}& 1\\ \cline{1-3}
  \multicolumn{1}{c}{2} & 1 & \multicolumn{1}{c}{1} & 4 
\end{array}
\hspace*{2mm}
\begin{array}{|@{}c@{}c@{}c@{}|c}\cline{1-3}
 x_{211} & x_{212}& x_{213}& 1\\
 x_{221} & x_{222}& x_{223}& 2\\
 x_{231} & x_{232}& x_{233}& 1\\ \cline{1-3}
  \multicolumn{1}{c}{1} & 2 & \multicolumn{1}{c}{1} & 4 
\end{array}
\hspace*{2mm}
\begin{array}{|@{}c@{}c@{}c@{}|c}\cline{1-3}
 x_{311} & x_{312}& x_{313}& 1\\
 x_{321} & x_{322}& x_{323}& 1\\
 x_{331} & x_{332}& x_{333}& 2\\ \cline{1-3}
  \multicolumn{1}{c}{1} & 1 & \multicolumn{1}{c}{2} & 4 
\end{array}
\hspace*{2mm}
\begin{array}{|@{}c@{}c@{}c@{}|c}\cline{1-3}
  2& 1& 1& 4\\
  1& 2& 1& 4\\
  1& 1& 2& 4\\ \cline{1-3}
  \multicolumn{1}{c}{4} & 4 & \multicolumn{1}{c}{4} & 12 
\end{array}
\]
where the rightmost table shows the marginal $\{x_{+jk}\}$. 
There are $18$ elements in this fiber as follows.
\[{\footnotesize
\hspace*{1mm}1: 
\begin{array}{|c|}\hline
2\ 0\ 0\\ 0\ 1\ 0\\ 0\ 0\ 1\\ \hline
\end{array}
\hspace*{1mm}
\begin{array}{|c|}\hline
0\ 1\ 0\\ 1\ 1\ 0\\ 0\ 0\ 1\\ \hline
\end{array}
\hspace*{1mm}
\begin{array}{|c|}\hline
0\ 0\ 1\\ 0\ 0\ 1\\ 1\ 1\ 0\\ \hline
\end{array}
\hspace*{0mm}
\hspace*{2mm}2: 
\begin{array}{|c|}\hline
2\ 0\ 0\\ 0\ 1\ 0\\ 0\ 0\ 1\\ \hline
\end{array}
\hspace*{1mm}
\begin{array}{|c|}\hline
0\ 1\ 0\\ 1\ 0\ 1\\ 0\ 1\ 0\\ \hline
\end{array}
\hspace*{1mm}
\begin{array}{|c|}\hline
0\ 0\ 1\\ 0\ 1\ 0\\ 1\ 0\ 1\\ \hline
\end{array}
\hspace*{0mm}
\hspace*{2mm}3: 
\begin{array}{|c|}\hline
2\ 0\ 0\\ 0\ 1\ 0\\ 0\ 0\ 1\\ \hline
\end{array}
\hspace*{1mm}
\begin{array}{|c|}\hline
0\ 1\ 0\\ 0\ 1\ 1\\ 1\ 0\ 0\\ \hline
\end{array}
\hspace*{1mm}
\begin{array}{|c|}\hline
0\ 0\ 1\\ 1\ 0\ 0\\ 0\ 1\ 1\\ \hline
\end{array}
}\]
\[{\footnotesize
\hspace*{0mm}4: 
\begin{array}{|c|}\hline
2\ 0\ 0\\ 0\ 1\ 0\\ 0\ 0\ 1\\ \hline
\end{array}
\hspace*{1mm}
\begin{array}{|c|}\hline
0\ 0\ 1\\ 1\ 1\ 0\\ 0\ 1\ 0\\ \hline
\end{array}
\hspace*{1mm}
\begin{array}{|c|}\hline
0\ 1\ 0\\ 0\ 0\ 1\\ 1\ 0\ 1\\ \hline
\end{array}
\hspace*{0mm}
\hspace*{2mm}5: 
\begin{array}{|c|}\hline
2\ 0\ 0\\ 0\ 0\ 1\\ 0\ 1\ 0\\ \hline
\end{array}
\hspace*{1mm}
\begin{array}{|c|}\hline
0\ 1\ 0\\ 1\ 1\ 0\\ 0\ 0\ 1\\ \hline
\end{array}
\hspace*{1mm}
\begin{array}{|c|}\hline
0\ 0\ 1\\ 0\ 1\ 0\\ 1\ 0\ 1\\ \hline
\end{array}
\hspace*{0mm}
\hspace*{2mm}6: 
\begin{array}{|c|}\hline
2\ 0\ 0\\ 0\ 0\ 1\\ 0\ 1\ 0\\ \hline
\end{array}
\hspace*{1mm}
\begin{array}{|c|}\hline
0\ 0\ 1\\ 0\ 2\ 0\\ 1\ 0\ 0\\ \hline
\end{array}
\hspace*{1mm}
\begin{array}{|c|}\hline
0\ 1\ 0\\ 1\ 0\ 0\\ 0\ 0\ 2\\ \hline
\end{array}
}\]
\[{\footnotesize
\hspace*{0mm}7: 
\begin{array}{|c|}\hline
1\ 1\ 0\\ 1\ 0\ 0\\ 0\ 0\ 1\\ \hline
\end{array}
\hspace*{1mm}
\begin{array}{|c|}\hline
1\ 0\ 0\\ 0\ 2\ 0\\ 0\ 0\ 1\\ \hline
\end{array}
\hspace*{1mm}
\begin{array}{|c|}\hline
0\ 0\ 1\\ 0\ 0\ 1\\ 1\ 1\ 0\\ \hline
\end{array}
\hspace*{0mm}
\hspace*{2mm}8: 
\begin{array}{|c|}\hline
1\ 1\ 0\\ 1\ 0\ 0\\ 0\ 0\ 1\\ \hline
\end{array}
\hspace*{1mm}
\begin{array}{|c|}\hline
1\ 0\ 0\\ 0\ 1\ 1\\ 0\ 1\ 0\\ \hline
\end{array}
\hspace*{1mm}
\begin{array}{|c|}\hline
0\ 0\ 1\\ 0\ 1\ 0\\ 1\ 0\ 1\\ \hline
\end{array}
\hspace*{0mm}
\hspace*{2mm}9: 
\begin{array}{|c|}\hline
1\ 1\ 0\\ 1\ 0\ 0\\ 0\ 0\ 1\\ \hline
\end{array}
\hspace*{1mm}
\begin{array}{|c|}\hline
0\ 0\ 1\\ 0\ 2\ 0\\ 1\ 0\ 0\\ \hline
\end{array}
\hspace*{1mm}
\begin{array}{|c|}\hline
1\ 0\ 0\\ 0\ 0\ 1\\ 0\ 1\ 1\\ \hline
\end{array}
}\]
\[{\footnotesize
10: 
\begin{array}{|c|}\hline
1\ 1\ 0\\ 0\ 0\ 1\\ 1\ 0\ 0\\ \hline
\end{array}
\hspace*{1mm}
\begin{array}{|c|}\hline
1\ 0\ 0\\ 0\ 2\ 0\\ 0\ 0\ 1\\ \hline
\end{array}
\hspace*{1mm}
\begin{array}{|c|}\hline
0\ 0\ 1\\ 1\ 0\ 0\\ 0\ 1\ 1\\ \hline
\end{array}
\hspace*{2mm}
11: 
\begin{array}{|c|}\hline
1\ 1\ 0\\ 0\ 0\ 1\\ 1\ 0\ 0\\ \hline
\end{array}
\hspace*{1mm}
\begin{array}{|c|}\hline
0\ 0\ 1\\ 1\ 1\ 0\\ 0\ 1\ 0\\ \hline
\end{array}
\hspace*{1mm}
\begin{array}{|c|}\hline
1\ 0\ 0\\ 0\ 1\ 0\\ 0\ 0\ 2\\ \hline
\end{array}
\hspace*{2mm}
12: 
\begin{array}{|c|}\hline
1\ 0\ 1\\ 1\ 0\ 0\\ 0\ 1\ 0\\ \hline
\end{array}
\hspace*{1mm}
\begin{array}{|c|}\hline
1\ 0\ 0\\ 0\ 2\ 0\\ 0\ 0\ 1\\ \hline
\end{array}
\hspace*{1mm}
\begin{array}{|c|}\hline
0\ 1\ 0\\ 0\ 0\ 1\\ 1\ 0\ 1\\ \hline
\end{array}
}\]
\[{\footnotesize
13: 
\begin{array}{|c|}\hline
1\ 0\ 1\\ 1\ 0\ 0\\ 0\ 1\ 0\\ \hline
\end{array}
\hspace*{1mm}
\begin{array}{|c|}\hline
0\ 1\ 0\\ 0\ 1\ 1\\ 1\ 0\ 0\\ \hline
\end{array}
\hspace*{1mm}
\begin{array}{|c|}\hline
1\ 0\ 0\\ 0\ 1\ 0\\ 0\ 0\ 2\\ \hline
\end{array}
\hspace*{2mm}
14: 
\begin{array}{|c|}\hline
1\ 0\ 1\\ 0\ 1\ 0\\ 1\ 0\ 0\\ \hline
\end{array}
\hspace*{1mm}
\begin{array}{|c|}\hline
1\ 0\ 0\\ 0\ 1\ 1\\ 0\ 1\ 0\\ \hline
\end{array}
\hspace*{1mm}
\begin{array}{|c|}\hline
0\ 1\ 0\\ 1\ 0\ 0\\ 0\ 0\ 2\\ \hline
\end{array}
\hspace*{2mm}
15: 
\begin{array}{|c|}\hline
1\ 0\ 1\\ 0\ 1\ 0\\ 1\ 0\ 0\\ \hline
\end{array}
\hspace*{1mm}
\begin{array}{|c|}\hline
0\ 1\ 0\\ 1\ 1\ 0\\ 0\ 0\ 1\\ \hline
\end{array}
\hspace*{1mm}
\begin{array}{|c|}\hline
1\ 0\ 0\\ 0\ 0\ 1\\ 0\ 1\ 1\\ \hline
\end{array}
}\]
\[{\footnotesize
16: 
\begin{array}{|c|}\hline
1\ 0\ 1\\ 0\ 1\ 0\\ 1\ 0\ 0\\ \hline
\end{array}
\hspace*{1mm}
\begin{array}{|c|}\hline
0\ 1\ 0\\ 1\ 0\ 1\\ 0\ 1\ 0\\ \hline
\end{array}
\hspace*{1mm}
\begin{array}{|c|}\hline
1\ 0\ 0\\ 0\ 1\ 0\\ 0\ 0\ 2\\ \hline
\end{array}
\hspace*{2mm}
17: 
\begin{array}{|c|}\hline
0\ 1\ 1\\ 1\ 0\ 0\\ 1\ 0\ 0\\ \hline
\end{array}
\hspace*{1mm}
\begin{array}{|c|}\hline
1\ 0\ 0\\ 0\ 2\ 0\\ 0\ 0\ 1\\ \hline
\end{array}
\hspace*{1mm}
\begin{array}{|c|}\hline
1\ 0\ 0\\ 0\ 0\ 1\\ 0\ 1\ 1\\ \hline
\end{array}
\hspace*{2mm}
18: 
\begin{array}{|c|}\hline
0\ 1\ 1\\ 1\ 0\ 0\\ 1\ 0\ 0\\ \hline
\end{array}
\hspace*{1mm}
\begin{array}{|c|}\hline
1\ 0\ 0\\ 0\ 1\ 1\\ 0\ 1\ 0\\ \hline
\end{array}
\hspace*{1mm}
\begin{array}{|c|}\hline
1\ 0\ 0\\ 0\ 1\ 0\\ 0\ 0\ 2\\ \hline
\end{array}
}\]
Now consider connecting these elements by the set of $9$ basic moves
 such as (\ref{eqn:basic-move-3way}). 
The undirected graph we obtain is Figure \ref{fig:graph-3x3x3-example}.
\begin{figure}[htbp]
\begin{center}
\includegraphics[width=80mm]{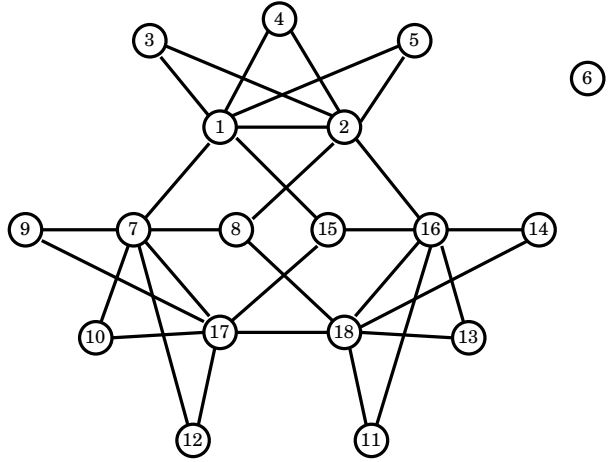}
\caption{Undirected graph obtained from the set of the basic moves.}
\label{fig:graph-3x3x3-example}
\end{center}
\end{figure}
Because this is not connected, the set of the basic moves is not a
Markov basis. This example shows that we need moves such as
\begin{equation}
 u_{111}u_{122}u_{133}u_{213}u_{221}u_{232} -
 u_{113}u_{121}u_{132}u_{211}u_{222}u_{233}
\label{eqn:deg6-moves}
\end{equation}
to constitute a Markov basis. 

Now calculate a Markov basis by Macaulay2 for this example. Using
 $a,b,c$ for the sufficient statistics instead of $v$, the following
is the commands to calculate a reduced Gr\"obner basis for this problem.
{\footnotesize
\begin{verbatim}
R = QQ[a11,a12,a13,a21,a22,a23,a31,a32,a33,
       b11,b12,b13,b21,b22,b23,b31,b32,b33,
       c11,c12,c13,c21,c22,c23,c31,c32,c33,
       x111,x112,x113,x121,x122,x123,x131,x132,x133,
       x211,x212,x213,x221,x222,x223,x231,x232,x233,
       x311,x312,x313,x321,x322,x323,x331,x332,x333,
       MonomialOrder=>{27,27}]
I = ideal(x111-a11*b11*c11,x112-a11*b12*c12,x113-a11*b13*c13,
          x121-a12*b11*c21,x122-a12*b12*c22,x123-a12*b13*c23,
          x131-a13*b11*c31,x132-a13*b12*c32,x133-a13*b13*c33,
          x211-a21*b21*c11,x212-a21*b22*c12,x213-a21*b23*c13,
          x221-a22*b21*c21,x222-a22*b22*c22,x223-a22*b23*c23,
          x231-a23*b21*c31,x232-a23*b22*c32,x233-a23*b23*c33,
          x311-a31*b31*c11,x312-a31*b32*c12,x313-a31*b33*c13,
          x321-a32*b31*c21,x322-a32*b32*c22,x323-a32*b33*c23,
          x331-a33*b31*c31,x332-a33*b32*c32,x333-a33*b33*c33)
G = gb(I); g = gens(G)
selectInSubring(1,g)
\end{verbatim}
}
Unfortunately, this calculation may be hard to carry out for average
 PC. In fact, I could not finish the above calculation within one hour by
 my slow laptop (with $2.80$ GHz CPU, $8.00$ GB RAM, running on
 vmware). Instead, check the calculation for $2\times 3\times 3$ cases. 
With the similar input commands, we have the output instantly in this
 case. 
From the output, we see that there are $1417$ elements in the 
reduced
 Gr\"obner basis of $I_A^*$, and $15$ elements in the 
reduced
 Gr\"obner basis of $I_A$ as follows.
{\footnotesize
\begin{verbatim}
i10 : selectInSubring(1,g)

o10 = | x122x133x223x232-x123x132x222x233 x112x133x213x232-x113x132x212x233
      -------------------------------------------------------------------------
      x121x133x223x231-x123x131x221x233 x121x132x222x231-x122x131x221x232
      -------------------------------------------------------------------------
      x111x133x213x231-x113x131x211x233 x111x132x212x231-x112x131x211x232
      -------------------------------------------------------------------------
      x112x123x213x222-x113x122x212x223 x111x123x213x221-x113x121x211x223
      -------------------------------------------------------------------------
      x111x122x212x221-x112x121x211x222
      -------------------------------------------------------------------------
      x112x121x133x211x223x232-x111x123x132x212x221x233
      -------------------------------------------------------------------------
      x111x122x133x213x221x232-x113x121x132x211x222x233
      -------------------------------------------------------------------------
      x111x122x133x212x223x231-x112x123x131x211x222x233
      -------------------------------------------------------------------------
      x113x121x132x212x223x231-x112x123x131x213x221x232
      -------------------------------------------------------------------------
      x112x121x133x213x222x231-x113x122x131x212x221x233
      -------------------------------------------------------------------------
      x111x123x132x213x222x231-x113x122x131x211x223x232 |

              1       15
o10 : Matrix R  <--- R
\end{verbatim}
}
We see that the set of the basic moves and the degree $6$ moves such as
 (\ref{eqn:deg6-moves}) actually constitutes a Markov basis for $2\times
 3\times 3$ cases. 

The calculation for $3\times 3\times 3$ cases can be carried out by
 faster software such as 4ti2 (\cite{4ti2}), that can be also used in
 Macaulay 2 as follows.
\begin{verbatim}
loadPackage "FourTiTwo"
A = matrix "1,1,1,0,0,0,0,0,0,0,0,0,0,0,0,0,0,0,0,0,0,0,0,0,0,0,0;
            0,0,0,1,1,1,0,0,0,0,0,0,0,0,0,0,0,0,0,0,0,0,0,0,0,0,0;
            0,0,0,0,0,0,1,1,1,0,0,0,0,0,0,0,0,0,0,0,0,0,0,0,0,0,0;
            0,0,0,0,0,0,0,0,0,1,1,1,0,0,0,0,0,0,0,0,0,0,0,0,0,0,0;
            0,0,0,0,0,0,0,0,0,0,0,0,1,1,1,0,0,0,0,0,0,0,0,0,0,0,0;
            0,0,0,0,0,0,0,0,0,0,0,0,0,0,0,1,1,1,0,0,0,0,0,0,0,0,0;
            0,0,0,0,0,0,0,0,0,0,0,0,0,0,0,0,0,0,1,1,1,0,0,0,0,0,0;
            0,0,0,0,0,0,0,0,0,0,0,0,0,0,0,0,0,0,0,0,0,1,1,1,0,0,0;
            0,0,0,0,0,0,0,0,0,0,0,0,0,0,0,0,0,0,0,0,0,0,0,0,1,1,1;
            1,0,0,1,0,0,1,0,0,0,0,0,0,0,0,0,0,0,0,0,0,0,0,0,0,0,0;
            0,1,0,0,1,0,0,1,0,0,0,0,0,0,0,0,0,0,0,0,0,0,0,0,0,0,0;
            0,0,1,0,0,1,0,0,1,0,0,0,0,0,0,0,0,0,0,0,0,0,0,0,0,0,0;
            0,0,0,0,0,0,0,0,0,1,0,0,1,0,0,1,0,0,0,0,0,0,0,0,0,0,0;
            0,0,0,0,0,0,0,0,0,0,1,0,0,1,0,0,1,0,0,0,0,0,0,0,0,0,0;
            0,0,0,0,0,0,0,0,0,0,0,1,0,0,1,0,0,1,0,0,0,0,0,0,0,0,0;
            0,0,0,0,0,0,0,0,0,0,0,0,0,0,0,0,0,0,1,0,0,1,0,0,1,0,0;
            0,0,0,0,0,0,0,0,0,0,0,0,0,0,0,0,0,0,0,1,0,0,1,0,0,1,0;
            0,0,0,0,0,0,0,0,0,0,0,0,0,0,0,0,0,0,0,0,1,0,0,1,0,0,1;
            1,0,0,0,0,0,0,0,0,1,0,0,0,0,0,0,0,0,1,0,0,0,0,0,0,0,0;
            0,1,0,0,0,0,0,0,0,0,1,0,0,0,0,0,0,0,0,1,0,0,0,0,0,0,0;
            0,0,1,0,0,0,0,0,0,0,0,1,0,0,0,0,0,0,0,0,1,0,0,0,0,0,0;
            0,0,0,1,0,0,0,0,0,0,0,0,1,0,0,0,0,0,0,0,0,1,0,0,0,0,0;
            0,0,0,0,1,0,0,0,0,0,0,0,0,1,0,0,0,0,0,0,0,0,1,0,0,0,0;
            0,0,0,0,0,1,0,0,0,0,0,0,0,0,1,0,0,0,0,0,0,0,0,1,0,0,0;
            0,0,0,0,0,0,1,0,0,0,0,0,0,0,0,1,0,0,0,0,0,0,0,0,1,0,0;
            0,0,0,0,0,0,0,1,0,0,0,0,0,0,0,0,1,0,0,0,0,0,0,0,0,1,0;
            0,0,0,0,0,0,0,0,1,0,0,0,0,0,0,0,0,1,0,0,0,0,0,0,0,0,1"
R = QQ[x111,x112,x113,x121,x122,x123,x131,x132,x133,
       x211,x212,x213,x221,x222,x223,x231,x232,x233,
       x311,x312,x313,x321,x322,x323,x331,x332,x333]
I = toricMarkov(A,R)
\end{verbatim}
This calculation is finished within 1 second by my laptop. 
From the output, we see that $27$ basic moves such as
 (\ref{eqn:basic-move-3way}) and $54$ moves of degree $6$ such as
(\ref{eqn:deg6-moves}) constitute a minimal Markov basis \footnote{The
 4ti2 command {\tt toricMarkov} gives a minimal Markov basis as the
 output. We can also obtain a Gr\"obner basis by the command
 {\tt toricGroebner}.}. Using this minimal Markov basis, we can
 construct a connected Markov chain for this fiber. The corresponding
 undirected graph is Figure \ref{fig:graph-3x3x3-example-2}.
\begin{figure}[htbp]
\begin{center}
\includegraphics[width=80mm]{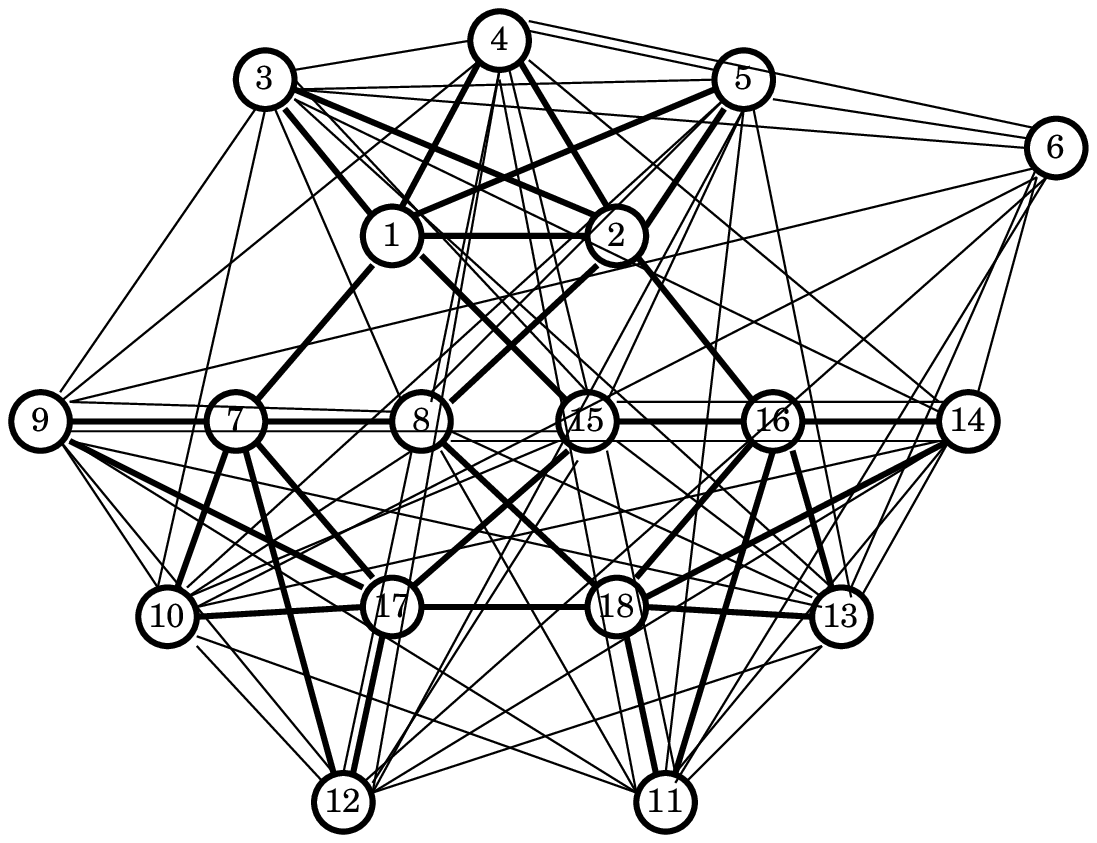}
\caption{Undirected graph obtained from a minimal Markov basis.}
\label{fig:graph-3x3x3-example-2}
\end{center}
\end{figure}
\end{Example}
Interestingly, for the problems of the larger sizes, the structure of
the Markov basis becomes more complicated. For example, for the no
three-factor interaction model of $3\times 3\times 4$ tables, the set of
degree $4, 6, 8$ moves becomes a Markov basis, and for 
$3\times 3\times 5$ tables, the set of 
degree $4, 6, 8, 10$ moves becomes a Markov basis. These results are
summarized in Chapter 9 of \cite{Aoki-Hara-Takemura-2012}.

\subsection*{Acknowledgments}
The author thanks the referee for his constructive suggestions that led
to significant improvement of this paper.

\bibliographystyle{plain}

\end{document}